%% file: AdmCoarsening.tex
\documentclass[3p]{elsarticle}

\input{general/packages}
\input{general/pgfplotsettings}

\input{general/commonCommands}

\drawboundingboxtrue
\drawboundingboxfalse 

\usepackage{amsmath,amsfonts}
\input{mathcommands}

\usepackage{svg}
\usepackage{gensymb}

\begin{document}

\begin{frontmatter}

\input{./general/titleAndAuthors}

%\input{./general/journalFooter} 

% Abstract ---------------------------------------
\vspace{-1.5cm} 
\hrule 

\input{./general/abstract}

%% Keywords ---------------------------------------
%\vspace{0.25cm}
%\noindent \textit{Keywords:} 
\input{general/keywords} 
%\vspace{0.25cm}
%\hrule 

\end{frontmatter}

%\tableofcontents

% Actual Content ---------------------------------
\input{./sections/introduction}
\input{./sections/ProblemFormulation.tex}
\input{./sections/THBSplines.tex}
\input{./sections/DiscretizedProblem.tex}
\input{./sections/Algorithms.tex}

\input{./sections/NumericalExample.tex}
\input{./sections/conclusions}
\section*{Acknowledgements} 
Massimo Carraturo and Alessandro Reali have been partially supported by Fondazione Cariplo - Regione Lombardia through the project ``Verso nuovi strumenti di simulazione super veloci ed accurati basati sull'analisi isogeometrica'', within the program RST - rafforzamento. Carlotta Giannelli has been partially supported by the Istituto Nazionale di Alta Matematica (INdAM) through Finanziamenti Premiali SUNRISE. Rafael V\'azquez has been partially supported by the {ERC} Advanced Grant ``CHANGE'', grant number 694515, 2016-2020.
Massimo Carraturo, Carlotta Giannelli, and Alessandro Reali are members of the INdAM Research group
GNCS (Gruppo Nazionale per il Calcolo Scientifico).  

%\input{general/acknowledgements} -------------------------------  

% Appendix
%\appendix
%\input{sections/datastructures} 

%% References -------------------------------------
\section*{\refname}
\bibliography{bibliography}

\end{document}

%% file: general/packages.tex
\usepackage{color} 
\usepackage{psfrag} 
\usepackage{subfig}
\usepackage{algorithm}
\usepackage{algorithmicx}
\usepackage{algpseudocode}
\usepackage{textcomp}
\usepackage{hyperref}
\usepackage{amsmath}
\usepackage{graphics}
 \usepackage{pstool} 
\usepackage[capitalise, noabbrev]{cleveref}
\usepackage[colorinlistoftodos]{todonotes}
\usepackage[export]{adjustbox} %/macros/latex/contrib/adjustbox

\usepackage{natbib}          
\setcitestyle{square,numbers,sort&compress}
\usepackage{pgfplots}
\usepackage{pgfplotstable}
\usepackage{filecontents}
\usepackage{pgf}

\usepackage{gensymb}
\usepackage[T1]{fontenc}
\usepackage{amsthm}
\usepackage[utf8]{inputenc}
\theoremstyle{definition}
\newtheorem{definition}{Definition}[section]
\theoremstyle{remark}

\usepackage{tikz}
\usetikzlibrary{decorations.pathmorphing}
\usepackage{tikzscale}
\usepackage{tikz-3dplot}
\usepackage{graphicx}

\usetikzlibrary{3d}
\usetikzlibrary{spy}
\usetikzlibrary{snakes,arrows,shapes}    
\usetikzlibrary{spy}
\usetikzlibrary{backgrounds}  
\usetikzlibrary{decorations}
\usetikzlibrary{positioning}
\usetikzlibrary{patterns}
\usetikzlibrary{calc} 
\usetikzlibrary{external} % Is activated with the command '\tikzexternalize'.
\usepackage{booktabs}   
\usepackage{makecell} 
\usepackage{colortbl}   
\usepackage{xspace}
\usepackage{color}     
\usepackage{setspace}

\usepackage{soul}
\usepackage{ulem}
\usepackage{currfile}
\usepackage{import}

\providecommand\phantomcaption{\caption@refstepcounter\@captype}

%% file: general/pgfplotsettings.tex
\pgfplotsset{compat=1.7}

\newcommand{\findmax}[3]{
    \pgfplotstablesort[sort key={#2},sort cmp={float >}]{\sorted}{#1}%
    \pgfplotstablegetelem{0}{#2}\of{\sorted}%
    \let #3=\pgfplotsretval%
}

\pgfplotscreateplotcyclelist{myCycleListBlackAndWhite}
{%
  {black, mark=o},
  {black, mark=square},
  {black, mark=triangle},
  {black, mark=diamond}, 
  {black, mark=pentagon}, 
  {black, mark=*},
  {black, mark=square*},
  {black, mark=triangle*},
  {black, mark=diamond*}, 
  {black, mark=pentagon*}, 
}

\definecolor{darkgreen}{rgb}{0,0.4,0} 
\definecolor{darkbrown}{rgb}{0.5, 0.396, 0.09}

\definecolor{c1}{rgb}{0.0, 0.4196078431372549, 0.6431372549019608}
\definecolor{c2}{rgb}{1.0, 0.5019607843137255, 0.054901960784313725}
\definecolor{c3}{rgb}{0.6705882352941176, 0.6705882352941176,
0.6705882352941176} \definecolor{c}{rgb}{0.34901960784313724, 0.34901960784313724, 0.34901960784313724}
\definecolor{c4}{rgb}{0.37254901960784315, 0.6196078431372549,
0.8196078431372549} \definecolor{c}{rgb}{0.7843137254901961, 0.3215686274509804, 0.0}
\definecolor{c5}{rgb}{0.5372549019607843, 0.5372549019607843,
0.5372549019607843} \definecolor{c}{rgb}{0.6352941176470588, 0.7843137254901961, 0.9254901960784314}
\definecolor{c6}{rgb}{1.0, 0.7372549019607844, 0.4745098039215686}
\definecolor{c7}{rgb}{0.8117647058823529, 0.8117647058823529,
0.8117647058823529}

\pgfplotscreateplotcyclelist{myCycleListColor}
{%
  {blue, mark=o},
  {red, mark=square},
  {darkbrown, mark=triangle},
  {darkgreen, mark=diamond},
  {black, mark=pentagon},
  {blue, mark=*},
  {red, mark=square*},
  {darkbrown, mark=triangle*},
  {darkgreen, mark=diamond*},
  {black, mark=pentagon*},
}

\pgfplotscreateplotcyclelist{myCycleListColorNoMarkers}
{
  {black},
  {blue, dashed},
  {red, dotted},
  {darkbrown, loosely dashed},
  {darkgreen, loosely dotted},
}

\pgfplotsset{every axis/.append style= 
              {
                font=\small,
                mark size=2,
                line width = 0.1,
                legend style={font=\small, mark size=3, draw=none, fill=none},
                legend cell align=left,
                cycle list name=myCycleListColor,
              }
            }
            
\pgfplotstableset
{
  every odd row/.style={before row={\rowcolor[gray]{0.9}}},
  every head row/.style=
  {
    before row=
    {%
      \toprule
    },
    after row=\midrule
  },
  every last row/.style=
  {
    after row=\bottomrule
  }
}

\newif\ifdrawboundingbox

\usetikzlibrary{external} % Is activated with the command
\tikzset{external/system call={pdflatex \tikzexternalcheckshellescape
-halt-on-error -interaction=batchmode -jobname "\image" "\texsource"}} 
\pgfkeys
{ 
  /pgf/number format/.cd,
  fixed,
  set thousands separator={\,},
  1000 sep in fractionals,
}

%% file: general/commonCommands.tex
\newcolumntype{C}[1]{>{\centering\arraybackslash}m{#1}}
\newcolumntype{R}[1]{>{\raggedright\arraybackslash}m{#1}}
\newcolumntype{L}[1]{>{\raggedleft\arraybackslash}m{#1}}

\newcommand{\delete}[1]{\xspace}

%% file: mathcommands.tex
\usepackage{etoolbox}
\ifdef{\R}{
	\renewcommand{\R}{\mathbb{R}}
}{
\newcommand{\R}{\mathbb{R}}
}
\ifdef{\N}{\renewcommand{\N}{\mathbb{N}}}{\newcommand{\N}{\mathbb{N}}}
\ifdef{\C}{\renewcommand{\C}{\mathbb{C}}}{\newcommand{\C}{\mathbb{C}}}
\ifdef{\Z}{\renewcommand{\Z}{\mathbb{Z}}}{\newcommand{\Z}{\mathbb{Z}}}

\newcommand{\lnorm}[1]{{\left\vert\kern-0.25ex\left\vert\kern-0.25ex\left\vert #1 
    \right\vert\kern-0.25ex\right\vert\kern-0.25ex\right\vert}}

\newcommand{\meinlemma}{Lemma \stepcounter{smeinTheorems}\the\value{smeinTheorems}}
\newcommand{\meinsatz}{Satz \stepcounter{smeinTheorems}\the\value{smeinTheorems}}
\newcommand{\meinedefinition}{Def.\ \stepcounter{smeinTheorems}\the\value{smeinTheorems}}

\newcounter{meineAussage} 
\newcommand{\opnorm}[1]{{\left\vert\kern-0.25ex\left\vert\kern-0.25ex\left\vert #1 
    \right\vert\kern-0.25ex\right\vert\kern-0.25ex\right\vert}}

\algnewcommand{\Oor}{\textbf{or}\,}
\algnewcommand{\Aand}{\textbf{and}\,}

%% file: general/titleAndAuthors.tex
\newcommand{\publicationDate}{\today}
\title{Suitably graded THB-spline refinement and coarsening: Towards an adaptive isogeometric analysis of additive manufacturing processes}

\author[PaviaAddress,MunichAddress]{Massimo Carraturo}
%\corref{mycorrespondingauthor}\cortext[mycorrespondingauthor]{Corresponding author}
  \ead{ massimo.carraturo01@universitadipavia.it}
\author[FirenzeAddress]{Carlotta Giannelli}
\ead{carlotta.giannelli@unifi.it}
\author[PaviaAddress]{Alessandro Reali}
\ead{alessandro.reali@unipv.it}
\author[LosannaAddress,ImatiAddress]{Rafael V{\'a}zquez}
\ead{rafael.vazquez@epfl.ch}
\address[PaviaAddress]{Department of Civil Engineering and Architecture, 
University of Pavia,
via Ferrata 3, 27100 Pavia, Italy}
\address[MunichAddress]{Chair of Computation in Engineering, 
Technische Universit{\"a}t M\"unchen,
Arcisstr. 21, 80333 M\"unchen, Germany} 
\address[FirenzeAddress]{Dipartimento di Matematica e Informatica "U.Dini", 
Universit\`a degli Studi di Firenze,
Viale Morgagni 67/A, 50134 Florence, Italy}
\address[LosannaAddress]{Institute of Mathematics, 
Ecole Polytechnique F{\'e}d{\'e}rale de Lausanne, Station 8, 1015 Lausanne, Switzerland}
\address[ImatiAddress]{Istituto di Matematica Applicata e Tecnologie Informatiche ``E. Magenes'' del CNR, via Ferrata 5, 27100, Pavia, Italy}
%  
%%\newcommand{\journal}{Some Fancy Journal}
%\newcommand{\publicationDate}{\today}

%% file: general/abstract.tex
\begin{abstract}
In the present work we introduce a complete set of algorithms to efficiently perform adaptive refinement and coarsening by exploiting truncated hierarchical B-splines (THB-splines) defined on suitably graded isogeometric meshes, that are called admissible mesh configurations.
We apply the proposed algorithms to two-dimensional linear heat transfer problems with localized moving heat source, as simplified models for additive manufacturing applications.
We first verify the accuracy of the admissible adaptive scheme with respect to an overkilled solution, for then comparing our results with similar schemes which consider different refinement and coarsening algorithms, with or without taking into account grading parameters. 
This study shows that the THB-spline admissible solution delivers an optimal discretization for what concerns not only the accuracy of the approximation, but also the (reduced) number of degrees of freedom per time step.
In the last example we investigate the capability of the algorithms to approximate the thermal history of the problem for a more complicated source path. 
The comparison with uniform and non-admissible hierarchical meshes demonstrates that also in this case our adaptive scheme returns the desired accuracy while strongly improving the computational efficiency. 
\end{abstract}

%% file: general/keywords.tex
\begin{keyword}
Isogeometric analysis \sep adaptivity \sep hierarchical splines \sep THB-splines \sep heat transfer analysis \sep additive manufacturing
\end{keyword}

%% file: sections/introduction.tex
\section{Introduction}
% IgA into
Isogeometric analysis (IGA)~\citep{hughes_isogeometric_2005, cottrell_isogeometric_2009} is a paradigm to solve partial differential equations employing smooth spline functions as basis for the analysis.
The classical  isoparametric approach of the finite element method is inverted and the exact geometry representation is now introduced within the analysis.
The original idea underlying this novel approach is to develop a tight connection between computer aided design (CAD) and numerical analysis, with the aim of performing analysis directly from CAD models to avoid tedious and time-consuming meshing processes. IGA has been so far successfully applied in many engineering fields including, among others, structural analysis~\citep{cottrell_studies_2007}, biomechanics~\citep{kamensky_immersogeometric_2015}, structural dynamics~\citep{Hughes_2014}, and contact mechanics~\citep{de_lorenzis_isogeometric_2014}.
\par

As CAD standard for spline representations, B-splines and non-uniform rational B-splines (NURBS) are the most commonly used spline technologies in the isogeometric setting.
Nevertheless, due to their tensor product structure, they are not well suited to treat localized phenomena.
% Local refinement in IGA
Hierarchical B-splines (HB-splines)~\citep{Forsey_1988,kraft1997} is an adaptive spline technology that enables the possibility to properly deal with local problems. Its application in isogeometric analysis has been widely studied, see e.g.~\citep{Vuong_2011,schillinger_isogeometric_2012,Scott_2014}.
Based on the multi-level concept of HB-splines, truncated hierarchical B-splines (THB-splines)~\citep{Giannelli_2012}, have also been proposed as an effective tool to perform hierarchical refinement while reducing the interactions between different levels in the spline hierarchy. The truncated basis has been successfully applied in different problems related to computer aided design~\citep{kiss2014b,bracco2018} and isogeometric analysis~\citep{giannelli2016,Henning_2016,dangella_multi-level_2018}.
It is worth to mention that, while several papers investigated refinement schemes for adaptive isogeometric methods in the last years, only very recently, few authors also focused on the study of suitable and effective mesh coarsening \cite{Lorenzo_2017,hennig2018}.
\par
% Additive manufacturing
Additive Manufacturing (AM) technologies are undergoing an exponential growth in many engineering fields, from aerospace to biomedical applications.
As a direct consequence of AM diffusion, there is an increasing demand for efficient and reliable numerical technologies not only to gain a deeper understanding of the physical process but also to optimize process parameters and predict the final shape of manufactured products.
\par
In particular, laser powder bed fusion (L-PBF) is a specific AM process where a laser beam selectively melts a layer of metal powder, building the final product by means of a layer-by-layer process.
In L-PBF technologies a single layer is~$\approx 30\mu\text{m}$ thick and the laser beam radius is around~$50\mu\text{m}$, while the dimensions of the product are in the order of centimeters.
This multi-scale nature in space makes the simulation of such a process an extremely challenging task.
Moreover, melting and solidification phenomena occurs in few microseconds in a very small region around the laser spot.
Considering that a complete L-PBF process lasts hours, we can easily understand that the computational costs for a high-fidelity simulation of such a process becomes too expensive even for the most recent supercomputers.
\par
Therefore, simulations of L-PBF processes represent a challenging task which involves many aspects of numerical, engineering and material science~\citep{king_laser_2015}.
% Open challanges in AM
In particular, we can identify three main open tasks in developing a validated and reliable simulation tool for this kind of processes: 
\begin{enumerate}
\item Develop a suitable physical and numerical model to represent the phenomena occurring during the process~\citep{lundback_modelling_2011};
\item Deliver accurate measurements of the physical parameters to validate the physical and numerical model~\citep{heigel_measurement_2017};
\item Develop a numerical technology which addresses the computational issues of the multi-physical and multi-scale nature of the problem~\citep{pal_generalized_2016,papadakis_computational_2014}. 
\end{enumerate} 
\par
% Scope of the work
The present work aims at introducing an efficient and adaptive numerical scheme based on THB-splines defined on suitably graded meshes to solve linear heat transfer problem including a localized heat source traveling on a two-dimensional domain. The THB-spline refinement routine follows the approach introduced in \cite{buffa_2016}, while the coarsening algorithm which automatically preserves the grading properties of the hierarchical mesh configuration is newly introduced. The results are obtained using GeoPDEs~\citep{de_falco_geopdes_2011,Vazquez_2016}, a MATLAB\textsuperscript \textregistered~ toolbox which implements, beside classical tensor product IGA, also adaptive hierarchical discretizations~\citep{bracco_2018,garau_2018}.
Even if our simple model neglects many features of the process, it still includes the spatial multi-scale issue present in L-PBF processes. 
The choice of employing an adaptive isogeometric scheme to solve the problem of a traveling laser source is justified by the presence of high thermal gradients in the region surrounding the laser beam. 
% Why IGA?
In~\citep{kollmannsberger_hierarchical_2018} it is demonstrated that high-order schemes combined together with adaptive mesh refinement and coarsening can optimally treat localized problems involving steep gradients, motivating the choice of employing an adaptive isogeometric method for this kind of applications.
\par
We outline that, even if we consider a simple, linear, two-dimensional model, the presented algorithms can be directly extended to 3D, non-linear, and multi-physics problems due to the extreme flexibility of the presented discretization technique.
Nevertheless, such an extension goes beyond the scope of the present work, whose goal is to show the effectiveness of adaptive isogeometric methods, based on THB-splines with suitable mesh grading, in the context of heat transfer problems with a localized heat source.
\par
% outline of the paper
The paper is organized as follows. In~\cref{sec:HeatTransferProblem} we present the equations for a linear heat transfer problem.
\cref{sec:IGA} briefly recalls the fundamental concepts of THB-splines, while~\cref{sec:HeatTransferDiscreteProblem} introduces the discrete form of the problem.
\cref{sec:AdapMeshAlg} presents the implementation of the method with a special focus on the coarsening procedure.
Then, in~\cref{sec:NumericalExamples}, we discuss two numerical examples with different scan paths and problem parameters, carrying out sensitivity studies for different levels of refinement and comparing results between uniform and adaptive discretizations, with or without considering suitable mesh grading.
Finally,~\cref{sec:conclusions} draws our concluding comments together with some final remarks.

%% file: sections/ProblemFormulation.tex
\section{Heat transfer problem} \label{sec:HeatTransferProblem}
In this section the strong and the weak form of the governing equations of a heat transfer problem are presented.
\par
Assuming the material obeys the Fourier's law of heat conduction with a Lagrangian reference frame, the problem can be written using the linear  heat transfer model as described in~\cite{bathe_finite_2007}.
Let us consider a temporal domain $T\subset\mathbb{R}$ and a spatial domain $\Omega\subset\mathbb{R}^{2}$ with Neumann boundaries $\Gamma_{N}$ such that $\Gamma_N=\partial\Omega$. The heat flow equilibrium in the interior of the body gives:
\begin{equation}
C_p\rho\dfrac{\partial\theta(\mathbf{x},t)}{\partial t} - \nabla\cdot(k\nabla\theta(\mathbf{x},t)) = f(\mathbf{x},t) \qquad \text{in}\quad\Omega\times T,
\label{eq:HeatTransfer}
\end{equation}
where $f$ is the external heat source, $k$ is the thermal conductivity, $C_p$ is the specific heat capacity and $\rho$ is the density of the solid material.
\par
The problem in equation~(\ref{eq:HeatTransfer}) is solved with respect to the temperature field $\theta=\theta(\mathbf{x},t)$ function of space and time, and defined under the following initial conditions
\begin{equation}
\theta(\mathbf{x},0) = \theta_{0}, \qquad \text{in}\quad\Omega,
\label{eq:initialCondition}
\end{equation}
and adiabatic boundary conditions
\begin{equation}
k\nabla\theta(\mathbf{x},t)\cdot\mathbf{n} = 0, \qquad \text{on}\quad\Gamma_{N}\times T,
\label{eq:NeumannBoundary}
\end{equation}
where $\theta_0$ is the initial temperature of the body and $\mathbf{n}$ is the exterior unit normal vector.
Note that we consider a simple linear heat transfer model in order to focus on the spatial discretization of the problem defined in equations~\eqref{eq:HeatTransfer}-\eqref{eq:NeumannBoundary}.
Consequently, we consider only the adiabatic boundary conditions defined in~\eqref{eq:NeumannBoundary} by neglecting radiation and convection effects on the domain boundaries.
\par
At time $t$ the weak solution of the problem in~(\ref{eq:HeatTransfer})-(\ref{eq:NeumannBoundary}) is obtained using the \textit{principle of virtual temperature} which can be written as follows:
\begin{equation}
C_p\rho\int_{\Omega}\tilde{\theta}(\mathbf{x},t)\dfrac{\partial\theta}{\partial t}(\mathbf{x},t)d\Omega + k\int_{\Omega}\nabla\tilde{\theta}(\mathbf{x},t)\cdot\nabla\theta(\mathbf{x},t)d\Omega = \int_{\Omega}\tilde{\theta}(\mathbf{x},t)f(\mathbf{x},t)d\Omega
\label{eq:principleOfVirtualTemperature}
\end{equation}
where $\tilde{\theta}$ is the virtual temperature.

%% file: sections/THBSplines.tex
\section{Truncated hierarchical B-splines} \label{sec:IGA}
This section presents the basic concepts of \textit{admissible} adaptive isogeometric methods, by exploiting truncated hierarchical B-splines \citep{Giannelli_2012} defined on suitably graded meshes as basis for the analysis. The definitions of this section mainly follow \cite{buffa_2016,buffa_2017b}, where a sound mathematical theory for adaptive isogeometric methods is fully developed, while related numerical results were recently presented in \cite{bracco2019}.
\subsection{Introduction to hierarchical spaces}
Given a nested sequence of parametric domains $\hat{\Omega}^0\supseteq\ldots\supseteq\hat{\Omega}^{N-1}$, subsets of a closed hyper-rectangle $D \in \mathbb{R}^d$, we can construct the hierarchical B-spline space of depth $N$ by considering a hierarchy of tensor-product B-spline bases $\mathcal{\hat{B}}^\ell$ of degree $p$ defined on the grid $\hat{G}^\ell$ for each level $\ell$, with $\ell=0,1,\ldots,N-1$. We assume that the domain $\Omega^\ell$ considered at level $\ell$ is the union of cells of the previous level $\ell-1$.
\begin{figure}[t]
\centerline{\includegraphics[width=.9\textwidth]{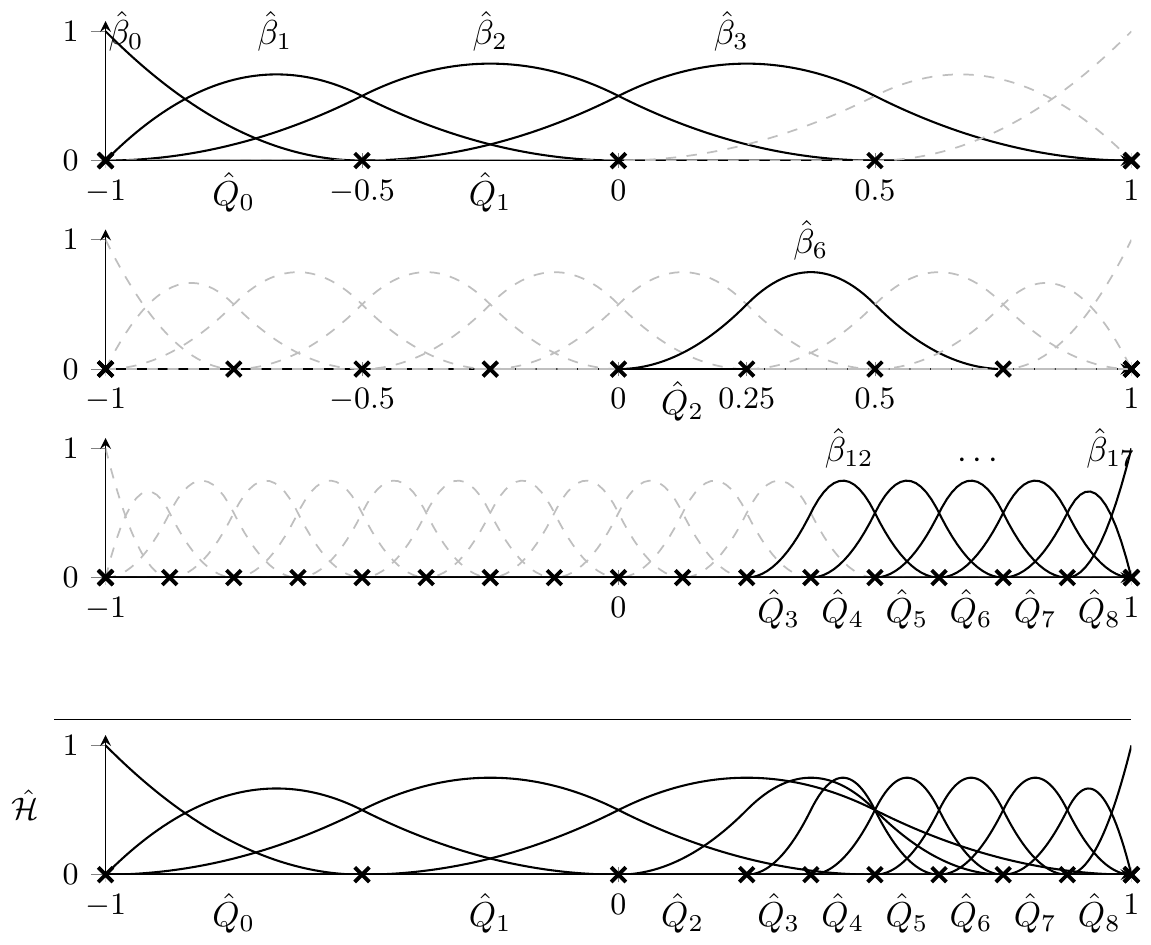}}
\caption{Univariate hierarchical B-splines defined on three refinement levels. The elements $\hat{Q}_i$ are selectively activated to refine the discretization towards the right-end of the domain. Active basis functions (solid lines) are shown together with inactive functions (dashed line).\label{fig:HierarchicalBasis}}
\end{figure}
A hierarchical mesh $\hat{\mathcal{Q}}$ collects the active cells which represent the elements of our discretization. It can be defined as
\begin{equation*}
\mathcal{\hat{Q}}:=
\left\lbrace \hat{Q}\in\hat{\mathcal{G}}^\ell, \ell=0,\ldots,N-1 \right\rbrace,
\end{equation*}
where $\hat{\mathcal{G}}^\ell$ is the set of active cells of level $\ell$, namely
\begin{equation*}
\mathcal{\hat{G}}^\ell:=\left\lbrace \hat{Q} \in \hat{G}^\ell:
\hat{Q} \subset \hat{\Omega}^\ell \wedge 
\hat{Q} \not\subset \hat{\Omega}^{\ell+1} \right\rbrace.
\end{equation*}
\cref{fig:HierarchicalBasis} shows an example of hierarchical mesh for $d=1$.~\cref{fig:MeshEvolution} and  \ref{fig:MultiTrackMeshEvolution} in~\cref{sec:NumericalExamples} instead illustrate several hierarchical mesh configurations in the bivariate setting. In the present work we consider only dyadic refinement, i.e., the children $\hat{Q_i}$, $i=1,\ldots,2^d$, of an active cell $\hat{Q}$ are obtained by bisection.
We finally define the hierarchical B-spline basis $\hat{\mathcal{H}}$ on the hierarchical mesh $\mathcal{\hat{Q}}$ as:
\begin{equation}
\hat{\mathcal{H}}(\hat{\mathcal{Q}}):=\left\lbrace \hat{\beta}\in\hat{\mathcal{B}}^\ell: \text{supp}\hat{\beta}\subseteq \hat{\Omega}^\ell \wedge \text{supp}\hat{\beta} \not \subseteq \hat{\Omega}^{\ell+1}, \ell=0,\ldots,N-1 \right\rbrace,
\label{eq:hierarchicalBasis}
\end{equation}
where supp$\hat{\beta}$ denotes the intersection of the support of $\hat{\beta}$ with $\hat{\Omega}^0$. Hierarchical B-splines for a univariate example are shown in \cref{fig:HierarchicalBasis}.
\subsection{The truncated basis}

\begin{figure}[t]
\centerline{\includegraphics[width=.9\textwidth]{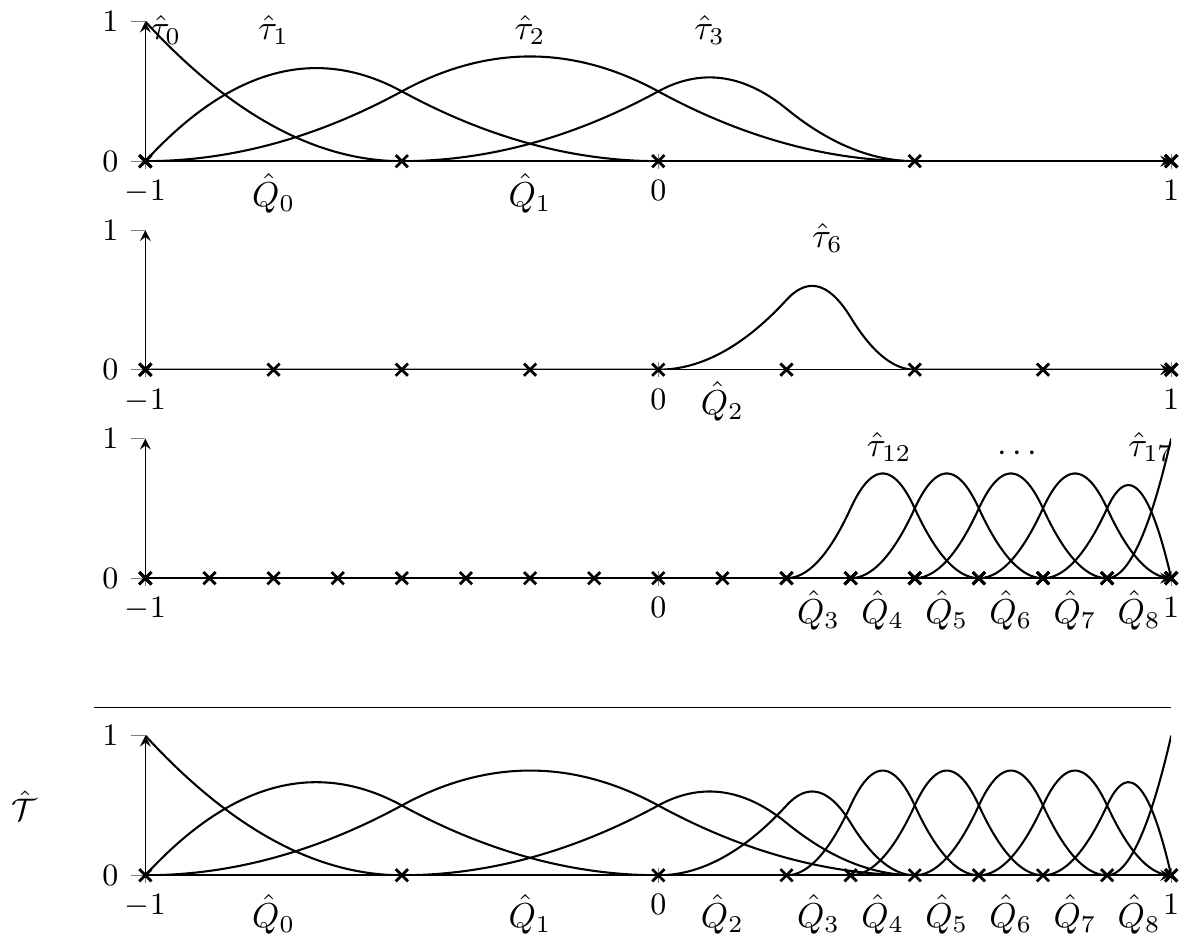}}
\caption{Univariate truncated hierarchical B-splines defined on three refinement levels. The elements $\hat{Q}_i$ are selectively activated to refine the discretization towards the right-end of the domain. Active truncated basis functions on each level  are shown. \label{fig:TruncatedHierarchicalBasis}}
\end{figure}

B-spline representations offer the possibility of suitably exploiting efficient refinement rules when nested spline spaces are considered. We can then consider the representation of $\hat{s}\in V^\ell\subset V^{\ell+1}$ with respect to the tensor-product B-spline basis $\hat{\mathcal{B}}^{\ell+1}$,
\[
\hat{s}=\sum_{\hat{\beta}\in \hat{\mathcal{B}}^{\ell+1}} c_{\hat{\beta}}^{\ell+1}(s)\hat{\beta},
\]
and define the truncation of $\hat{s}$ with respect to level ${\ell+1}$ as follows:
\begin{equation*}
\text{trunc}^{\ell+1}\hat{s}:=\sum_{\hat{\beta}\in \hat{\mathcal{B}}^{\ell+1},\text{supp}\hat{\beta}\not\subseteq\hat{\Omega}^{\ell+1}} c_{\hat{\beta}}^{\ell+1}(s)\hat{\beta},
\end{equation*} 
where $c_{\hat{\beta}}^{\ell+1}(s)$ is the coefficient of the function $\hat{s}$ with respect to the basis function $\hat{\beta}$ at level $\ell+1$.
By iteratively applying the truncation operation at the hierarchical basis functions in $\hat{\mathcal{H}}$, we obtain the truncated hierarchical basis (see~\cref{fig:TruncatedHierarchicalBasis}).
\par
\definition{
\textit{ 
The truncated hierarchical B-spline (THB-spline) basis $\hat{\mathcal{T}}$ with respect to the mesh $\hat{\mathcal{Q}}$ is defined as
\begin{equation*}
\hat{\mathcal{T}}(\hat{\mathcal{Q}}):=
\left\lbrace
\textnormal{Trunc}^{\ell+1}\hat{\beta}:\hat{\beta}\in \hat{\mathcal{B}}^\ell\cap\hat{\mathcal{H}}(\hat{\mathcal{Q}}), \ell=0,\ldots,N-2
\right\rbrace
\bigcup
\left\{
\hat{\beta}^{N-1}: \hat{\beta}^{N-1}\in \hat{\mathcal{B}}^{N-1}\cap \hat{\mathcal{H}}(\hat{\mathcal{Q}})
\right\}
,
\end{equation*}
where $\textnormal{Trunc}^{\ell+1}\hat{\beta}:=\textnormal{trunc}^{N-1}
(\textnormal{trunc}^{N-2}
(\ldots(\textnormal{trunc}^{\ell+1}(\hat{\beta} ) )\ldots ) )$,
for any $\hat{\beta}\in\hat{\mathcal{B}}^\ell \cap \hat{\mathcal{H}}(\hat{\mathcal{Q}})$.
}
}
\par
The key basic properties of the truncated basis are the following: non-negativity, linear independence, partition of unity, and, in addition, $\textnormal{span}\,\hat{\mathcal{H}}=\textnormal{span}\,\hat{\mathcal{T}}$ \citep{Giannelli_2012}.
\subsection{Admissible meshes}\label{sec:admissible_meshes}
In this work we use the definition of \textit{classes of admissible meshes} introduced in~\citep{buffa_2016}.
\par
\definition{\textit{A mesh $\hat{\mathcal{Q}}$ is admissible of class m if the truncated basis functions in $\hat{\mathcal{T}}(\hat{\mathcal{Q}})$ which take non-zero values over any element $\hat{Q}\in\hat{\mathcal{Q}}$ belong to at most m successive levels. }}
\par
This definition allows to consider hierarchical meshes where the number of THB-splines acting on a single element does not depend on the total number of levels in the hierarchy, but only on the parameter $m$.
\par
In order to implement the algorithms for admissible refinement and coarsening, see~\cref{sec:RefAlg} and~\cref{sec:CoarseAlg}, we need to consider three additional sets of elements: the \textit{multilevel support extension} of an element, together with the so-called \textit{refinement and coarsening neighborhoods}.

\definition{\textit{The multilevel support extension of an element $\hat{Q}\in \hat{G}^\ell$ with respect to level $k$, with $0\leq k \leq \ell$, is defined as:}}
\begin{equation*}
S(\hat{Q},k):=\left\lbrace \hat{Q'}\in\hat{G}^k:\exists\hat{\beta}\in \mathcal{\hat{B}}^k, \text{supp}\hat{\beta} \cap \hat{Q'} \neq \emptyset \wedge \text{supp}\hat{\beta} \cap \hat{Q} \neq \emptyset \right\rbrace.
\end{equation*}
That is, the multilevel support extension is formed by the support of B-splines of level $k$, such that they do not vanish on an element $\hat{Q}$.
By following \citep{buffa_2016}, we now consider the neighborhood of an element for the refinement algorithm, which we  rename as \textit{refinement neighborhood} to differentiate it from the analogous set used in coarsening.

\definition{\textit{The refinement neighborhood $\mathcal{N}_r(\mathcal{\hat{Q}},\hat{Q},m)$ of an element $\hat{Q}$ of level $\ell$ with respect to the class of admissibility $m$ is  defined as:}}
\begin{equation*}
\mathcal{N}_r(\mathcal{\hat{Q}},\hat{Q},m):=\left\lbrace \hat{Q}'\in\mathcal{\hat{G}}^{\ell-m+1}:\exists\, \hat{Q}''\in S(\hat{Q}, \ell-m+2), \hat{Q}''\subseteq \hat{Q}' \right\rbrace.
\end{equation*}
Note that, when considering the element $\hat{Q}$ of level $\ell$ to be refined into four elements of level $\ell+1$, we want to guarantee that any THB-spline of level $\ell-m+1$ acting on $\hat{Q}$ is fully truncated with respect to level $\ell-m+2$, so that it will vanish on the children of $\hat{Q}$ which will be activated after the refinement of $\hat{Q}$. By (recursively) refining the elements in the multilevel support extension of level $\ell-m+2$
, as in the definition of the refinement neighborhood, guarantees that the admissibility property is maintained. Note however that this choice is conservative in nature, since we do not consider the real support of the THB-splines, but only their truncation with respect to level $\ell-m+2$, while they could be truncated at intermediate levels depending on the mesh configuration.
%
%Note that, by considering the multilevel support extension with respect to level $\ell-m+2$ in the definition of the refinement neighborhood, we can check the smallest set of elements that should be recursively refined when the element $\hat{Q}$ of level $\ell$ is subdivided into four elements of level $\ell+1$ if we want to guarantee that any THB-spline of level $\ell-m+1$ acting on $\hat{Q}$ is completely truncated with respect to level $\ell-m+2$. These THB-splines will then be zero on the children on $\hat{Q}$ in the refined mesh and admissibility property is preserved. Note however that this choice is conservative in nature since any THB-spline could also be truncated at intermediate levels in different mesh configurations. 

Finally, we introduce here the definition of \textit{coarsening neighborhood} of an element, an additional set of elements required for admissible coarsening.
\definition{\textit{The coarsening neighborhood $\mathcal{N}_c(\mathcal{\hat{Q}},\hat{Q},m)$ of an element $\hat{Q}$ of level $\ell$ with respect to the class of admissibility $m$ is defined as:}}
\begin{equation*}
\mathcal{N}_c(\mathcal{\hat{Q}},\hat{Q},m):=\left\lbrace \hat{Q}'\in\mathcal{\hat{G}}^{\ell+m}:\exists\, \hat{Q}'\subset S(\hat{Q}'', \ell+1),\textnormal{ with } \hat{Q}''\in\mathcal{\hat{G}}^{\ell+1}\textnormal{ and } \hat{Q}''\subset \hat{Q} \right\rbrace.
\end{equation*}
The set $\mathcal{N}_c(\mathcal{\hat{Q}},\hat{Q},m)$ includes the active elements of level $\ell + m$ that are contained in the support extension of the children of $\hat{Q}$. By starting from an admissible mesh of class $m$, when considering the element $\hat{Q}$ of level $\ell$ to be reactivated, to preserve the admissibility condition we have to ensure that newly activated basis functions will vanish on the elements of level $\ell+m$. This is guaranteed when the coarsening neighborhood is empty, because newly added functions will be fully truncated with respect to level $\ell+1$. This property allows us to define a fully automatic way to preserve admissibility of the mesh, that we exploit in the algorithm of Section~\ref{sec:CoarseAlg}.
%The set $\mathcal{N}_c(\mathcal{\hat{Q}},\hat{Q},m)$ includes the elements of level $\ell + m$ that are contained in the support extension of the children of $\hat{Q}$. By starting from an admissible mesh of class $m$ and avoiding to reactivate an element whose coarsening neighborhood is empty we can easily define a fully automatic way to preserve the admissibility of the mesh with respect to the considered value of $m$.

%% file: sections/DiscretizedProblem.tex
\section{Adaptive isogeometric analysis} \label{sec:HeatTransferDiscreteProblem}
In this section the discrete form of the problem of~\cref{sec:HeatTransferProblem} is derived.
%---------------------------------------------------------------------------------------------------------------------------
\subsection{Spatial discretization} \label{sec:FEM}
In order to define the mesh and basis functions for the physical domain $\Omega$ we introduce the isogeometric mapping $\mathbf{F}:\hat{\Omega}^0\rightarrow\Omega$, such that:
\begin{equation*}
\mathbf{x}\in\Omega,
\quad \mathbf{x}=\mathbf{F}(\hat{\mathbf{x}})=\sum_{\hat{\tau}\in\hat{\mathcal{T}}_0}\mathbf{C}_{\hat{\tau}}\hat{\tau}\left(\hat{\mathbf{x}}\right),
\end{equation*}
with $\hat{\mathbf{x}}\in\hat{\Omega}^0$, $\mathbf{C}_{\hat{\tau}}\in\mathbb{R}^d$, and where $\hat{\mathcal T}_0$ is the truncated basis defined on an initial tensor-product mesh.
\par
The corresponding hierarchical mesh in the physical domain can be written as:
\begin{equation*}
\mathcal{Q}=\left\lbrace Q=\mathbf{F}(\hat{Q}):\hat{Q}\in\hat{\mathcal{Q}} \right\rbrace
\end{equation*}
and, analogously,
\begin{equation*}
\Omega^\ell=\mathbf{F}(\hat{\Omega}^\ell),
\quad \mathcal{G}^\ell=\left\lbrace Q\in\mathcal{Q}:\hat{Q}\in\mathcal{\hat{G}}^\ell \right\rbrace\quad 
\text{and} 
\quad G^\ell=\left\lbrace Q\subset\Omega:\hat{Q}\in\hat{G}^\ell \right\rbrace.
\end{equation*}
Equation \eqref{eq:principleOfVirtualTemperature} can now be written in a discrete form applying the isogeometric expansion to the temperature field $\theta$, such that:
\begin{equation*}
\theta(\mathbf{x},t) = \mathbf{N}(\mathbf{x})\boldsymbol{\theta}_t,
\end{equation*}
where $\boldsymbol{\theta}_t$ is the column vector of temperature degrees of freedom (DOFs) of the hierarchical mesh $\mathcal{Q}$ at time $t$ and $\mathbf{N}(\mathbf{x})$ is the row vector of the corresponding THB-spline basis functions. 
We can now rewrite the integrals of equation~\eqref{eq:principleOfVirtualTemperature} at a given instant in time in a matrix form, such that:
\begin{equation}
\mathbf{M}\boldsymbol{\dot\theta}_t + \mathbf{K}\boldsymbol{\theta}_t = \mathbf{f}_t,
\label{eq:FEMVirtualTemperature}
\end{equation}
where 
\begin{equation*}
\mathbf{M}=C_p\rho\int_{\Omega}\mathbf{N}^T\mathbf{N}d\Omega,\quad \mathbf{K}=\kappa\int_{\Omega}\mathbf{B}^T\mathbf{B}d\Omega,
\end{equation*}
and
\begin{equation*}
\mathbf{f}_t=\int_{\Omega}\mathbf{N}^Tf(\mathbf{x},t)d\Omega,
\end{equation*}
with $\mathbf{B}(\mathbf{x})=\nabla\mathbf{N}(\mathbf{x})$.
Equation~\eqref{eq:FEMVirtualTemperature} presents the time derivative of the temperature vector, therefore an additional discretization in time is required to numerically solve the problem.
%---------------------------------------------------------------------------------------------------------------------------
\subsection{Time integration} \label{sec:TimeDisc}
The linear system in~\eqref{eq:FEMVirtualTemperature} is approximated in time using the unconditionally stable backward Euler approach. 
At time $t + \Delta t$ equation~\eqref{eq:FEMVirtualTemperature} can be written in terms of the temperature increment $\Delta\boldsymbol{\theta}_{t+\Delta t}$ as:
\begin{equation}
\mathbf{M}\Delta\boldsymbol{\theta}_{t+\Delta t} + {\Delta t}\mathbf{K}\Delta\boldsymbol{\theta}_{t+\Delta t} = {\Delta t}\mathbf{f}_{t+\Delta t} - {\Delta t}\mathbf{K}\boldsymbol{\theta}_{t},
\label{eq:BackwardEuler}
\end{equation}
where $\Delta t$ is the constant time increment at each time step and $\boldsymbol{\theta}_t$ is the the solution vector at the previous time step.
Finally, we can iteratively solve the discrete problem defined in equation~\eqref{eq:BackwardEuler} using the algorithms described in the next section.
%-----------------------------------------------------------------
\subsection{Error estimator} \label{ssec:ErrorEstimator}
Before describing the implementation of the numerical method, we introduce the error estimator used to drive our adaptive scheme. We employ a residual-based \textit{a posteriori} error estimator $\varepsilon_{\mathcal{Q}}$ defined as:
\begin{equation}
\varepsilon_{\mathcal{Q}}^2(\theta,\mathcal{Q})= \sum_{Q\in\mathcal{Q}} \varepsilon^2_Q ,
\label{eq:errorEst}
\end{equation}
where 
\begin{equation*}
\quad\varepsilon^2_Q = h_Q^2 \int_Q \mid f_{t+\Delta t} -
C_p\rho\dfrac{\theta_{t+\Delta t}-\theta_t}{\Delta t} + \nabla\cdot(k\nabla\theta_{t+\Delta t})\mid^2 d Q  
\label{eq:errorEst1}
\end{equation*}
and $h_Q$ is the size of the element $Q$.
\par
Finally, to estimate the quality of the results with respect to a reference solution $\boldsymbol{\theta}_{t,ref}$, we need to define the error with respect to a certain norm. 
To this end we employ two different definitions of the energy of the system. 
Following~\citep{ohara_2011}, we define the internal energy of the system $\mathcal{E}_i(\boldsymbol{\theta}_t,\boldsymbol{\theta}_t)$ at time $t$ as
\begin{equation*}
\mathcal{E}_i(\boldsymbol{\theta}_t,\boldsymbol{\theta}_t)=
\dfrac{k}{2}\int_{\Omega}\nabla\theta_{t}\cdot\nabla\theta_{t}d\Omega,
\label{eq:internalEnergy}
\end{equation*}
whereas the total energy of the system $\mathcal{E}_T(\boldsymbol{\theta}_t,\boldsymbol{\theta}_t)$ at time $t$ is defined as
\begin{equation*}
\mathcal{E}_T(\boldsymbol{\theta}_t,\boldsymbol{\theta}_t)=
\dfrac{1}{2}\left(k\int_{\Omega}\nabla\theta_{t}\cdot\nabla\theta_{t}d\Omega +
C_p\rho\int_{\Omega}\theta_{t}\dfrac{\partial\theta_{t}}{\partial t}d\Omega \right).
\label{eq:totalEnergy}
\end{equation*}
We can then define the relative error in internal and total energy at a given instant in time $t$ as
\begin{equation*}
\varepsilon_{i}=
\sqrt{\dfrac
{\mid \mathcal{E}_i(\boldsymbol{\theta}_{t,ref},\boldsymbol{\theta}_{t,ref})-{\mathcal{E}_i(\boldsymbol{\theta}_{t},\boldsymbol{\theta}_{t})\mid^2}}
{\mid \mathcal{E}_i(\boldsymbol{\theta}_{t,ref},\boldsymbol{\theta}_{t,ref})\mid^2}
}
%\label{eq:internalEnergyError}
%\end{equation}
%%
\qquad\text{and}\qquad
%%
%\begin{equation}
\varepsilon_{T}=
\sqrt{\dfrac
{\mid \mathcal{E}_T(\boldsymbol{\theta}_{t,ref},\boldsymbol{\theta}_{t,ref})-{\mathcal{E}_T(\boldsymbol{\theta}_{t},\boldsymbol{\theta}_{t})\mid^2}}
{\mid \mathcal{E}_T(\boldsymbol{\theta}_{t,ref},\boldsymbol{\theta}_{t,ref})\mid^2}
},
%\label{eq:totalEnergyError}
\end{equation*}
respectively.

%% file: sections/Algorithms.tex
\section{Algorithms for admissible adaptivity} \label{sec:AdapMeshAlg}
This section presents the implementation of admissible adaptivity that we developed in GeoPDEs starting from the algorithms previously implemented in the code and described in~\citep{garau_2018,bracco_2018}.
The initial set of algorithms in~\citep{garau_2018} is now modified, since we aim here at solving the parabolic problem of equation~\eqref{eq:BackwardEuler} by employing admissible adaptive discretizations.
%---------------------------------------------------------------------------------------------------------------------------------------------------------------------------------------------------------
\subsection{Admissible adaptive backward Euler} \label{sec:AdapAlg}
Given an initial tensor-product mesh $\mathcal{Q}^0_0$, the corresponding truncated hierarchical B-spline space $\mathcal{T}_0^0$, an initial solution vector $\boldsymbol{\theta}_0^0$, a class of admissibility $m$, a maximum number of refinement levels $N$, and a tolerance \mbox{tol}, Algorithm~\ref{AdpMshAlg} returns an approximated solution of the problem defined in~\cref{sec:HeatTransferProblem} employing a backward Euler time integration scheme together with adaptive mesh refinement and coarsening procedures fulfilling the admissibility requirements. 
In this way, the algorithm allows to concentrate the computational efforts where the estimated errors are higher and, at the same time, to obtain an admissible mesh that gives the possibility to avoid undesired oscillations in the solution. 
The algorithm can be split, at each time step, into two separate parts: first, we refine the mesh iteratively, and, subsequently, we coarsen the mesh to generate the initial mesh for the next time step.
\floatname{algorithm}{Algorithm}
\renewcommand{\algorithmicrequire}{\textbf{Input:}}
\renewcommand{\algorithmicensure}{\textbf{Output:}}
\begin{algorithm}
\caption{\texttt{solve\_heat\_transfer\_problem}
%: Backward Euler time integration scheme using an adaptive mesh refinemement and coarsening fulfilling the admissibility requirements.
}\label{AdpMshAlg}
\begin{algorithmic}[1]
\Require{$\mathcal{Q}_0$, $\mathcal{T}_0$, $\boldsymbol{\theta}_0$, $m$, $N$, \mbox{tol}}
\Ensure{$\boldsymbol{\Theta}:=\left\lbrace\boldsymbol{\theta}_0,\ldots,\boldsymbol{\theta}_{t_{end}}\right\rbrace$}
\State{$t\gets0$} %\Comment{Initialize time}

\State{$\left(\mathcal{Q}_{t+\Delta t},\mathcal{T}_{t+\Delta t},  \varepsilon_{\mathcal{Q}}, \boldsymbol{\theta}_{t+\Delta t} \right)\gets$\texttt{adpt\_iter\_refine}($\mathcal{Q}_t$, $\mathcal{T}_t$, $\boldsymbol{\theta}_t$, $m$, $t$, $N$, \mbox{tol})}\Comment{~\cref{AdpIterRefAlg} ($1^{st}$ time step)}

\State{$t\gets t+\Delta t$}%\Comment{update time step}

\While{$t < t_{end}$}

\State{$\left(\mathcal{Q}_{t+\Delta t},\mathcal{T}_{t+\Delta t},  \varepsilon_{\mathcal{Q}}, \boldsymbol{\theta}_{t+\Delta t} \right)\gets$\texttt{adpt\_iter\_refine}($\mathcal{Q}_t$, $\mathcal{T}_t$, $\boldsymbol{\theta}_t$, $m$, $t$, 2, \mbox{tol})}\Comment{~\cref{AdpIterRefAlg}}

\State{$\mathcal{M}_c\gets$\texttt{mark\_min}($\varepsilon_{\mathcal{Q}}$,$\mathcal{Q}_{t+\Delta t}$)\Comment{~\cref{eq:MARK_COARSE}}}

\State{$\mathcal{Q}_{t+\Delta t}\gets$ \texttt{coarsen}($\mathcal{Q}_{t+\Delta t}$,$\mathcal{M}_c$,$m$)\Comment{~\cref{Coarse}}}

\State{$(\mathcal{T}_{t+\Delta t},\boldsymbol{\theta}_{t+\Delta t})\gets$ \texttt{project}($\mathcal{Q}_{t+\Delta t}$,$\mathcal{T}_{t+\Delta t},\boldsymbol{\theta}_{t+\Delta t}$)\Comment{$L_2$ projection onto the coarsened mesh}}

\State{$t\gets t + \Delta t$}

\EndWhile

\end{algorithmic}
\end{algorithm}
%---------------------------------------------------------------------------------------------------------------------------------------------------------------------------------------------------------
\subsection{Refinement} \label{sec:RefAlg}
\begin{algorithm}
\caption{\texttt{adpt\_iter\_refine}
%: Adaptively and iteratively refine the mesh and the corresponding function space. It returns the refined mesh and function space together with the corresponding solution vector and the estimated error.
}\label{AdpIterRefAlg}
\begin{algorithmic}[1]
\Require{$\mathcal{Q}^0$, $\mathcal{T}^0$, $\boldsymbol{\theta}_t$, $m$, $t$, $I_{\text{MAX}}$}, \mbox{tol}
\Ensure{$\left( \mathcal{Q}, \mathcal{T}, \varepsilon_{\mathcal{Q}}, \boldsymbol{\theta}_{t+\Delta t} \right)$}
\State{$i\gets0$} %\Comment{Initialize adaptivity iteration counter}

\State{$\Delta\boldsymbol{\theta}_{t+\Delta t}\gets$\texttt{solve}$(\mathcal{Q}^i,\mathcal{T}^i,\boldsymbol{\theta}_t)$}\Comment{~\cref{eq:BackwardEuler}}
	
\State{$\boldsymbol{\theta}_{t+\Delta t}\gets\boldsymbol{\theta}_{t}+\Delta\boldsymbol{\theta}_{t+\Delta t}$} %\Comment{Update solution}
	
\State{$\varepsilon_{\mathcal{Q}}\gets$\texttt{estimate}$(\boldsymbol{\theta}_{t+\Delta t},\mathcal{Q}^i,\mathcal{T}^i)$}\Comment{~\cref{eq:errorEst}}

\While{$i\leq I_{\text{MAX}}$ $\Aand$ $\varepsilon_{\mathcal{Q}} \geq \mbox{tol}$ }

	    \State{$\mathcal{M}_r\gets$\texttt{mark\_max}($\varepsilon_{\mathcal{Q}}$,$\mathcal{Q}^i$)}\Comment{~\cref{eq:MARK_REF}}
		
		\State{$\mathcal{Q}^{i+1}\gets$ \texttt{refine}($\mathcal{Q}^i$,$\mathcal{M}_r$,$m$)}\Comment{~\cref{REFINE}}

		\State{($\mathcal{T}^{i+1},\boldsymbol{\theta}_{t})\gets$ \texttt{project}($\mathcal{Q}^{i+1}$,$\mathcal{T}^i,\boldsymbol{\theta}_{t}$)\Comment{knot insertion as described in~\citep[Sec.~4.3]{garau_2018}}}
		
		\State{$i\gets i+1$} %\Comment{Update adaptivity iteration counter}
		
		\State{$\Delta\boldsymbol{\theta}_{t+\Delta t}\gets$\texttt{solve}$(\mathcal{Q}^i,\mathcal{T}^i,\boldsymbol{\theta}_t)$}
	
		\State{$\boldsymbol{\theta}_{t+\Delta t}\gets\boldsymbol{\theta}_{t}+\Delta\boldsymbol{\theta}_{t+\Delta t}$}
	
		\State{$\varepsilon_{\mathcal{Q}}\gets$\texttt{estimate}$(\boldsymbol{\theta}_{t+\Delta t},\mathcal{Q}^i,\mathcal{T}^i)$}
		
\EndWhile
\State{$\mathcal{Q}\gets\mathcal{Q}^{i}$, $\mathcal{T}\gets\mathcal{T}^{i}$}
\end{algorithmic}
\end{algorithm}
Algorithm~\ref{AdpMshAlg} at each time step calls the function~\texttt{adpt\_iter\_refine} defined in~\cref{AdpIterRefAlg} to adaptively refine the mesh and the corresponding function space.
The \texttt{adpt\_iter\_refine} function returns a refined mesh and the corresponding function space together with the solution vector and the estimated error on the refined space.
\cref{AdpIterRefAlg} starts solving equation~\eqref{eq:BackwardEuler} on the mesh obtained in the previous time step $\mathcal{Q}^0$. Successively, the algorithm estimates the element error and marks a set of active elements to be refined
\begin{equation}
\mathcal{M}_r=\texttt{mark\_max}\left(\varepsilon_{\mathcal{Q}}(\theta_t,Q)_{Q\in\mathcal{Q}}, \mathcal{Q}\right),
\label{eq:MARK_REF}
\end{equation}
following the so-called \textit{D\"orfler marking},
i.e. by considering a fixed refinement marking parameter $\alpha_r\in\left(0,1\right]$ such that
\begin{equation*}
\varepsilon_{\mathcal{Q}}({\theta_t,\mathcal{M}_r})\geq\alpha_r\varepsilon_{\mathcal{Q}}(\theta_t,\mathcal{Q}),
\end{equation*}
where $\theta_t$ is the discrete solution at time $t$.
\par%
The marked elements are then refined employing~\cref{REFINE} and~\cref{REFINE_RECURSIVE}, previously introduced in~\citep{buffa_2016}, see also \citep{bracco_2018} for a more detailed explanation of these algorithms which generate admissible meshes. Note that the choice of the parameter $m$ naturally influences the grading of the the hierarchical meshes.
\\
\begin{minipage}{\textwidth}
\begin{minipage}{.45\textwidth}
\begin{algorithm}[H]
\captionof{algorithm}{\texttt{refine}
%: given a set of active elements to be refined $\mathcal{M}_r$, update a mesh $\mathcal{Q}$ to get an admissible mesh $\mathcal{Q}^{\star}$ of class $m$.
}\label{REFINE}
\begin{algorithmic}[1]
\Require{$\mathcal{Q}$, $\mathcal{M}_r$, $m$}
\Ensure{$\mathcal{Q}^{\star}$}
\For{${Q}\in\mathcal{Q}\cap\mathcal{M}_r$}
\State{$\mathcal{Q}\gets$\texttt{refine\_recursive}($\mathcal{Q}$,$Q$,$m$)}	
\EndFor
\State{$\mathcal{Q}^{\star}\gets\mathcal{Q}$} %update solution
\end{algorithmic}
\end{algorithm}
\end{minipage}
\hfill%-------------------------------------------------------
\begin{minipage}{.45\textwidth}
\begin{algorithm}[H]
\captionof{algorithm}{\texttt{refine\_recursive}
%: recursively update the mesh $\mathcal{Q}$ to obtain an admissible refined mesh.
}\label{REFINE_RECURSIVE}
\begin{algorithmic}[1]
\Require{$\mathcal{Q}$, $Q$, $m$}
\Ensure{$\mathcal{Q}$}
\For{${Q'}\in\mathcal{N}_r(\mathcal{Q},Q,m)$}
\State{$\mathcal{Q}\gets$\texttt{refine\_recursive}($\mathcal{Q}$,$Q'$,$m$)}	
\EndFor
\If{$Q$ has not been subdivided}
\State{subdivide $Q$ and}
\State{update $\mathcal{Q}$ by replacing Q with its children}
\EndIf
\end{algorithmic}
\end{algorithm}
\end{minipage}
\end{minipage}
\bigskip

\cref{AdpIterRefAlg} is broken when either a maximum number of iterations is reached or the estimated error $\varepsilon_{\mathcal{Q}}$ is below a certain tolerance $\mbox{tol}$.
In particular, in this work we iterate the first time step until the maximum level of refinement $N$ is achieved, while, for the remaining time steps, we set the maximum number of iterations equal to two, see the call to \texttt{adpt\_iter\_refine} (\cref{AdpIterRefAlg}) on line 5 of \texttt{solve\_heat\_transfer\_problem} (\cref{AdpMshAlg}).
This choice is justified by the small values of $\Delta t$ that are generally used in AM simulations to discretize the problem in time.
Even if for linear problems a different choice for the time step size and, then, for the stopping criteria could be made, we aim at developing and verify an algorithm suitable to be used in more complex AM applications.
In these applications, in fact, the time step increments are necessarily very small due to the strong non linearity of the problem, therefore two consecutive meshes would not differ too much from each other.
\par
In view of the above considerations, at each iteration of~\cref{AdpIterRefAlg} the solution obtained at the previous time step on the coarse mesh is written in terms of basis functions of the locally refined mesh (using the knot insertion as described in~\citep{garau_2018}). A more accurate (and also more expensive) solution would consist instead in performing a global $L_2$ projection of the refined temperature vector of the previous time step, i.e. $\boldsymbol{\theta}_t$ before coarsening.
% we also choose to locally \textcolor{blue}{refine} (using the knot insertion as described in~\citep{garau_2018}) at each iteration of~\cref{AdpIterRefAlg} the solution obtained at the previous time step on the coarse mesh \textcolor{red}{usare la knot insertion non \`e una proiezione, potete riformulare la frase?}. 
%
%A more conservative (and also more expensive) solution would consist instead in performing a global $L_2$ projection of the refined temperature vector of the previous time step, i.e. $\boldsymbol{\theta}_t$ before coarsening.
%
In fact, with our choice we lose some accuracy each time we coarsen the mesh, but, since we employ very small time steps, this approximation remains acceptable as it will be shown in~\cref{sec:NumericalExamples}.
%---------------------------------------------------------------------------------------------------------------------------------------------------------------------------------------------------------
\input{./sections/coarsening}

%% file: sections/coarsening.tex
\subsection{Coarsening} \label{sec:CoarseAlg}
In order to complete the time step routine we need to coarsen the hierarchical mesh, that will be used in the next time step.
Also in this case we first mark a set of active elements to be coarsened
\begin{equation}
\mathcal{M}_c=\texttt{mark\_min}\left(\varepsilon_{\mathcal{Q}}(\theta_t,Q)_{Q\in\mathcal{Q}}, \mathcal{Q}\right),
\label{eq:MARK_COARSE}
\end{equation}
fixing a coarsening marking parameter $\alpha_c\in\left(0,1\right]$ and considering the elements with the lowest estimated error, such that
\begin{equation*}
\varepsilon_{\mathcal{Q}}({\theta_t,\mathcal{M}_c})\leq\alpha_c\varepsilon_{\mathcal{Q}}(\theta_t,\mathcal{Q}).
\end{equation*}
The last step is the coarsening of the mesh, which is described in Algorithm~\ref{Coarse}. An important issue is to decide, from the list of marked elements, the elements that should be reactivated. In principle, the parents of all marked elements could be reactivated. However, since coarsening implies a loss of information we have chosen to be conservative, and an element is reactivated only if all its children are marked. In other words, elements that are not marked will remain active. This is ensured by the first condition in line~6 of the algorithm.
Moreover, to guarantee that the coarsened mesh fulfills the admissibility property, we perform one more check on the selected elements: the element can be reactivated only if the coarsening neighborhood is empty, otherwise the admissibility condition would be violated, as we have explained in Section~\ref{sec:admissible_meshes}. If the chosen element satisfies the two conditions, the last step, performed in line~7 of Algorithm~\ref{Coarse}, updates the mesh by reactivating the element and removing its children.
\begin{algorithm}
\caption{\texttt{coarsen}
%: given a mesh $\mathcal{Q}$ and a set of active marked elements $\mathcal{M}_{c}$ to be coarsened, it returns an admissible coarsened mesh of class $m$.
}\label{Coarse}
\begin{algorithmic}[1]
\Require{$\mathcal{Q}$, $\mathcal{M}_c$, $m$}
\Ensure{$\mathcal{Q}$}
\For{$Q\in\mathcal{M}_c$}
\State $\mathcal{R}_{c} \gets \mathcal{R}_{c} \cup$ \texttt{get\_parent}($Q$)
\EndFor

\For{$Q\in\mathcal{R}_{c}$} \Comment{This loop must be done from the finest to the coarsest level}
\State $Q_c \gets$ \texttt{get\_children} ($Q$)
\If {($Q_c \subset \mathcal{M}_c$ \Aand $\mathcal{N}_c(\mathcal{Q},Q,m) = \emptyset$)}
\State{update $\mathcal{Q}$ by activating $Q$ and removing its children $Q_c$}
%\State $\mathcal{R}_{c,adm} \gets \mathcal{R}_{c,adm} \setminus Q$
%%\State{SE$\gets$\texttt{get\_support\_extension}($Q_c$,$\mathcal{Q}$)}
%%\State{Descendants$\gets$\texttt{get\_descendants}(SE,$m$)} \Comment{And the descendants have level $l+m$}
%%\If{Descendants $\cap$ $\mathcal{G}^{l+m}\neq\emptyset$}
%%{\color{red} \If{any Descendant is active or deactivated}} \Comment{Alternative to previous line}
\EndIf
\EndFor
\end{algorithmic}
\end{algorithm}

We remark that the coarsening algorithm for THB-splines that we propose differs from other algorithms presented in previous papers. In particular, in \cite{hennig2018} a global {\textit{a posteriori}} check for the admissibility of the mesh is performed after coarsening, which may cause to refine again reactivated elements. In our algorithm this check is local and performed before reactivating the elements. This guarantees an \textit{a priori} automatic control which saves computational time. In \cite{Lorenzo_2017} the admissibility condition is replaced by a similar concept, that the authors call \textit{function support balancing}, that also guarantees a certain grading of the mesh. As we will see in the numerical tests of \cref{sec:NumericalExamples}, the condition of admissibility maintains the same accuracy as the function support balancing, while the obtained refinement is more local for the former, which leads to spaces with many less degrees of freedom.

%% file: sections/NumericalExample.tex
\section{Numerical examples} \label{sec:NumericalExamples}
Two numerical examples are discussed in the following section.
In the first example, defined by a Gaussian heat source tavelling on a circular arc scan track, we compare the solution obtained using an admissible hierarchical mesh with the ones computed employing  uniform refinement, non-admissible hierarchical refinement and the function support balancing algorithm introduced in~\citet{Lorenzo_2017}.
%
%This comparison is attained evaluating the accuracy of the method for diffrent levels of refinement with respect to an \textit{overkilled} reference solution. 
%
The second example consists of a Gaussian heat source traveling on a multi-track source path. In this case we investigate the ability of the presented admissible adaptive scheme to capture the influence of adjacent scan tracks on the temperature evolution. Again we compare the solution of the admissible adaptive mesh with respect to the results obtained with uniform and non-admissible adaptive refinement.
In both the examples we use standard Gauss integration to numerically approximate the integrals of~\cref{eq:FEMVirtualTemperature} and we set the admissibility parameter $m$ equal to 2.
All the simulations in this work are performed using MATLAB\textsuperscript \textregistered on an Intel\textsuperscript \textregistered ~Core\textsuperscript \texttrademark ~i7-6700, CPU@3.40GHz, RAM 24Gb.
%---------------------------------------------------------------------------------------------------------------------------
\subsection{Circular arc scan track} \label{sec:CircularArc}
\begin{figure}[th]
\centering
\includegraphics[width=0.45\textwidth]{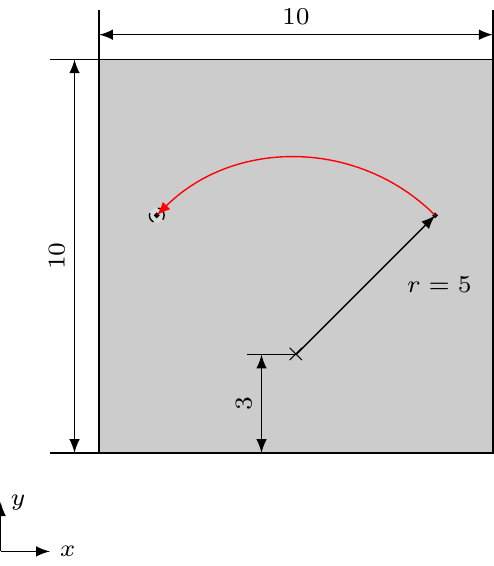}
\caption{Circular arc scan track: Heat source moving on a circular arc, the red arrow indicates the scan path (all distances are in mm).\label{ArcPathSetup}}
\end{figure}

The first example consists of a Gaussian heat source traveling along a circular arc on a $10\times10$ mm$^2$ surface domain, as described in~\cref{ArcPathSetup}. 
The external heat source function $f$ is modeled using a Gaussian function defined as:
\begin{equation}
f=\text{P}\eta\exp\left(-\left((x-x_0)^2+(y-y_0)^2\right)/r_h^2\right),
\label{rhsThermal}
\end{equation}
where P is the heat source power,
$\eta$ the absorptivity of the material, 
$r_h$ the heat source spot radius, 
while $(x_0,y_0)$ is the position of the source at a given instant in time.
In~\cref{ParameterTable} the values adopted for the process and the material model parameters are reported. In this first example we aim at comparing different discretization techniques without focusing on the real physics of the process, justifying in such a way the choice of simple unitary values for the material parameters.
\begin{table}[h]
\centering
\begin{tabular}{ | c | c | } 
\hline
\textbf{Parameters} & \textbf{Values} \\ 
\hline
Laser power P & $9\times 10^5$[W] \\ 
\hline
Laser speed & $1.57\left[\text{mm}/\text{sec}\right]$ \\ 
\hline
Absorptivity $\eta$ & $0.33$ \\ 
\hline
Source radius $r_h$ & $100\left[\mu\text{m}\right]$ \\
\hline
Conductivity $k$ & $1.0\left[\text{W}/\text{mm}/\text{K}\right]$ \\
\hline
Specific heat capacity $C_p$ & $1.0\left[\text{J}/\text{kg}/\text{K}\right]$ \\
\hline
Density $\rho$ & $1.0\left[\text{kg}/\text{mm}^3\right]$ \\
\hline
Initial temperature $\theta_0$ & $20.0\left[^{\circ}\text{C}\right]$ \\

\hline
\end{tabular}
\caption{Circular arc scan track: Process and material parameters.\label{ParameterTable}}
\end{table}
\par
For such a problem we generate an overkilled reference solution obtained with B-spline of degree 4 defined on $2^9\times2^9$ isogeometric tensor-product elements, while adaptive IGA discretization starts from a single knot span and is recursively bisected towards the regions with the highest error, as indicated by the error estimator defined in~\eqref{eq:errorEst} choosing $\alpha_r=0.1$ and $\alpha_c=0.25$. 
\par
\begin{figure}[t!]
\centering
\includegraphics[width=0.45\textwidth]{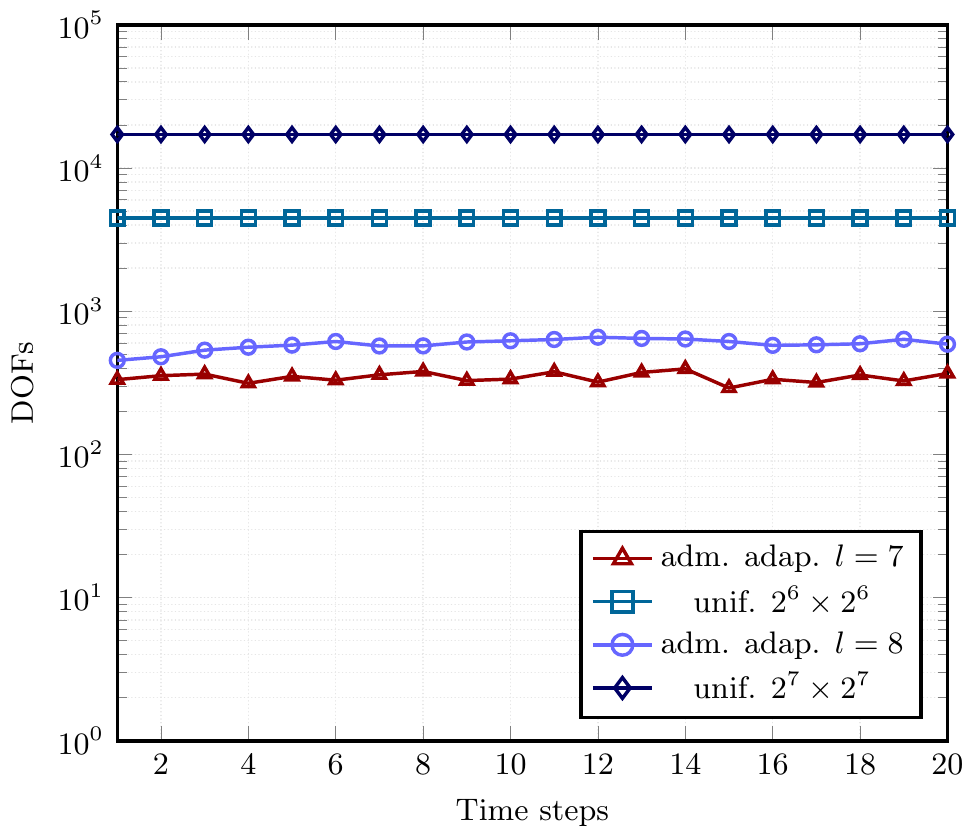}
\hfill
\includegraphics[width=0.45\textwidth]{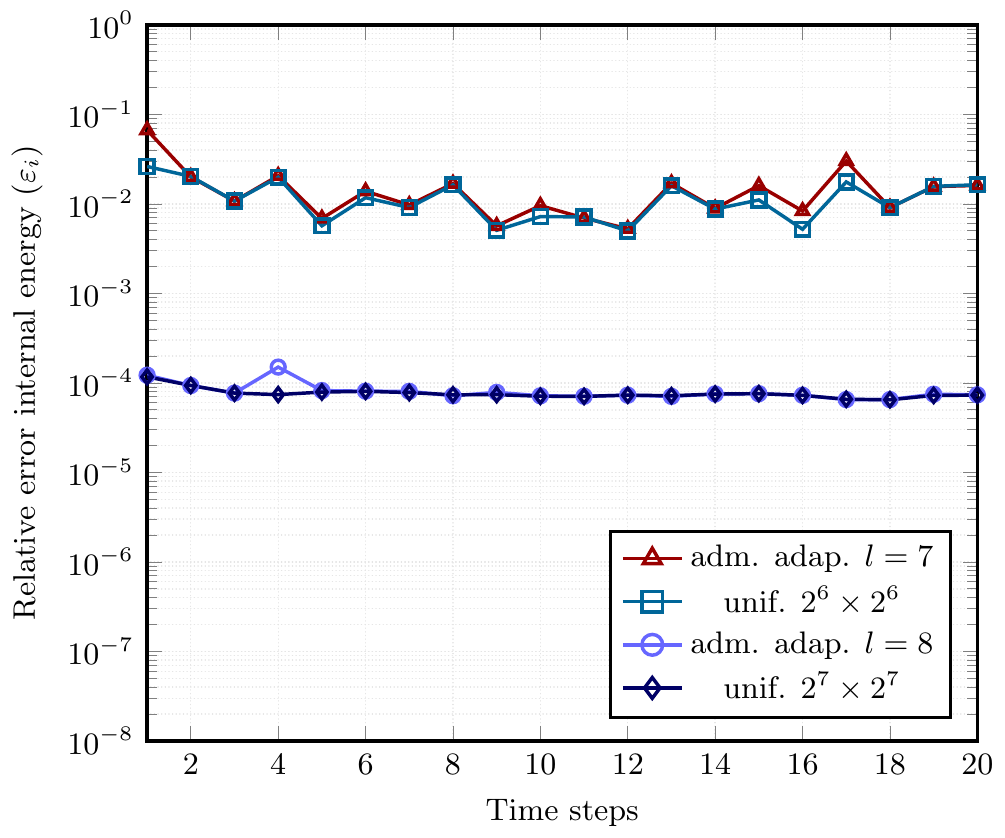}
\caption{Circular arc scan track: DOFs (right) and relative errors (right) comparison at each time step for different levels of refinement using adaptive THB-splines and uniform tensor-product meshes.}
\label{Exm1DOFS}
\end{figure}

\cref{Exm1DOFS} (right) shows the relative error in the internal energy norm $\varepsilon_i$ with respect to the reference solution for two different cubic admissible adaptive discretizations with seven and eight levels of refinement (adm. adap. $l=7$ and adm. adap. $l=8$) and their corresponding uniform meshes (unif. $2^6\times2^6$ and unif. $2^7\times2^7$), i.e. uniform meshes with elements of the same size of the smallest element in the adaptive mesh. 
It can be observed that the relative errors of the adaptive and the uniform discretization are almost identical for both cases, whereas, as shown in~\cref{Exm1DOFS} (left), the adaptive schemes require almost two orders of magnitude less DOFs compared to uniform IGA meshes. 
We want to stress the fact that this difference increases together with the maximum refinement level, i.e. the more localized is the problem the smaller is the resulting linear system (and thus the memory consumption) compared to the uniform case for the same level of accuracy.
% 
%All these requirements are addressed by the presented adaptive IGA implementation.  
%
\begin{figure}[t!]
\centering
%\subfloat[\label{DOFsCompAdm}]{% 
\includegraphics[width=0.45\textwidth]{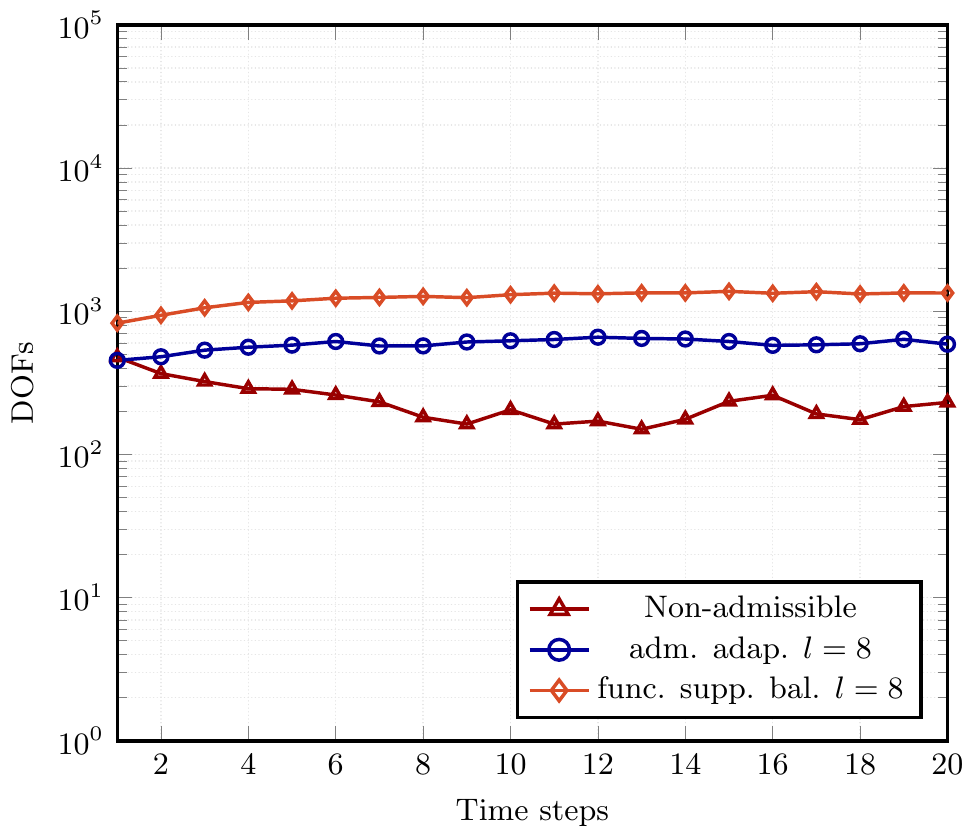}
  %}
% 
\hfill 
%\subfloat[\label{ErrorCompAdm}]{% 
\includegraphics[width=0.45\textwidth]{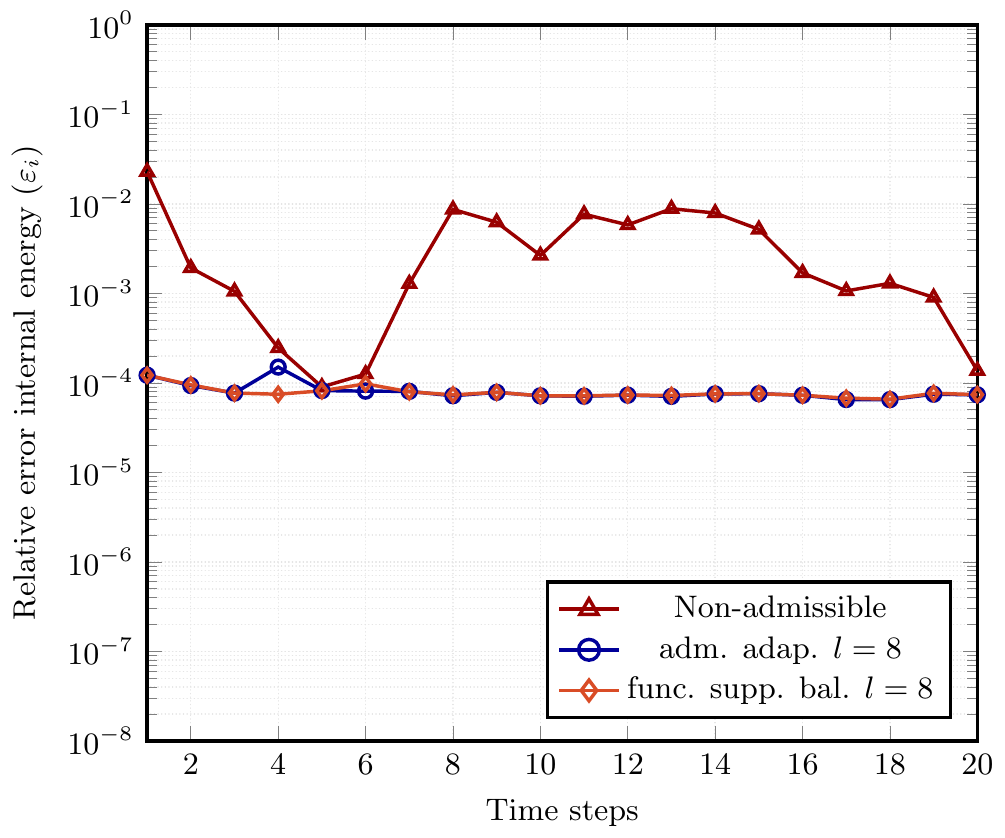}
%} 
 \caption{Circular arc scan track: DOFs (left) and relative errors (right) for the admissible adaptive mesh, the function support balancing~\citep{Lorenzo_2017} and a non-admissible mesh.}\label{NotAdmComparisonPlots}
\end{figure}
\begin{figure}[h!]
\centering
%\subfloat[\label{Temp1TS}]{% 
\includegraphics[width=0.329\textwidth, trim= 20mm 20mm 80mm 20mm,clip=true]{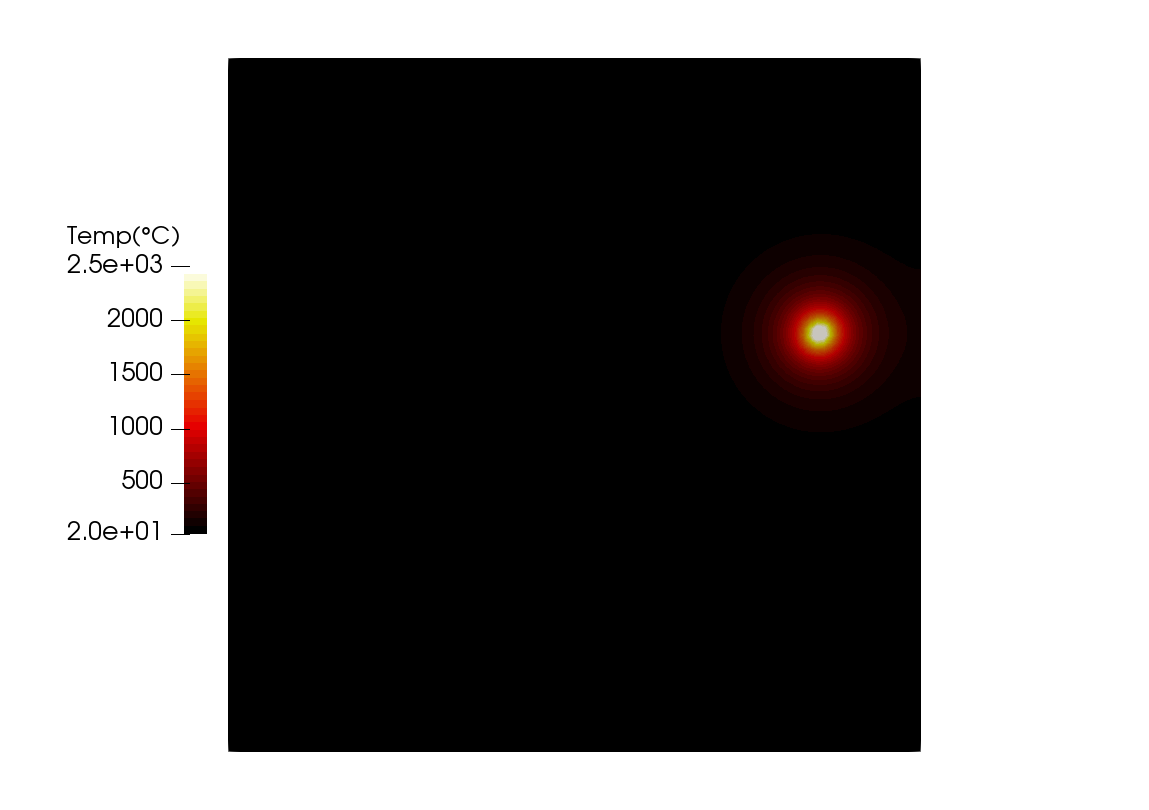}
%} 
%\subfloat[\label{Temp10TS}]{% 
\includegraphics[width=0.329\textwidth, trim= 20mm 20mm 80mm 20mm,clip=true]{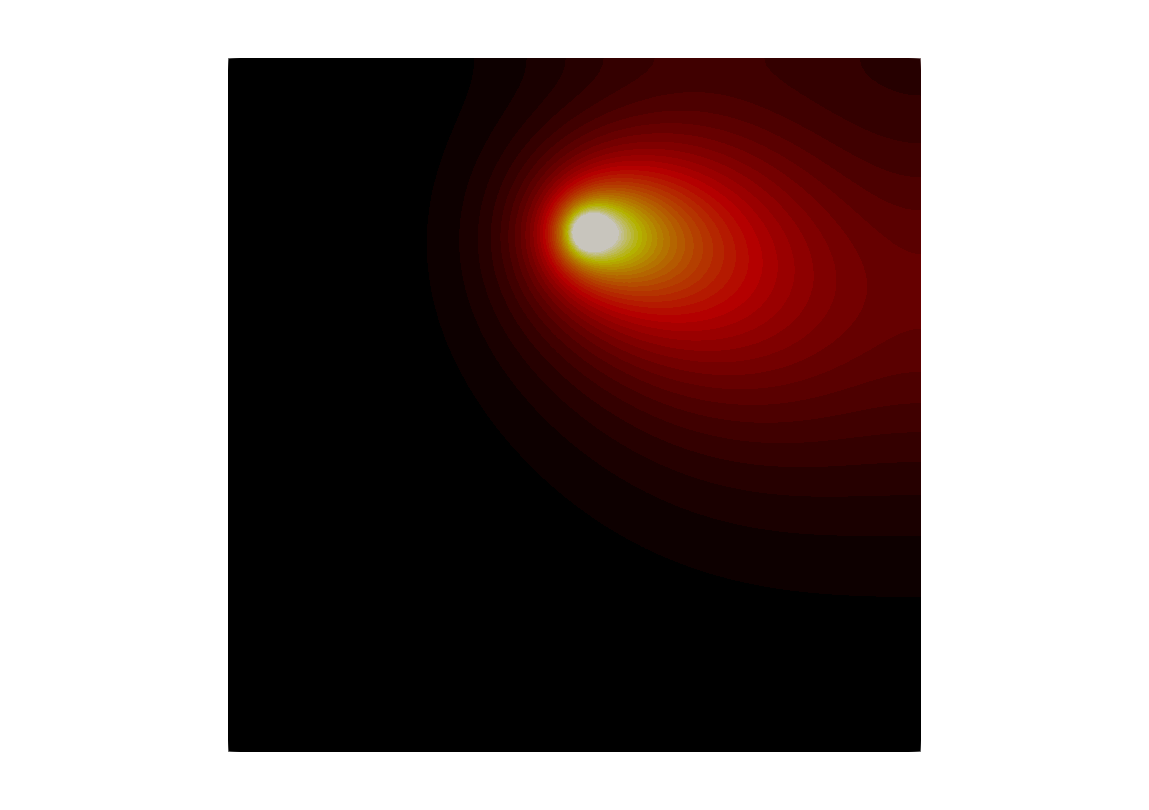}
%}
%\subfloat[\label{Temp20TS}]{% 
\includegraphics[width=0.329\textwidth, trim= 20mm 20mm 80mm 20mm,clip=true]{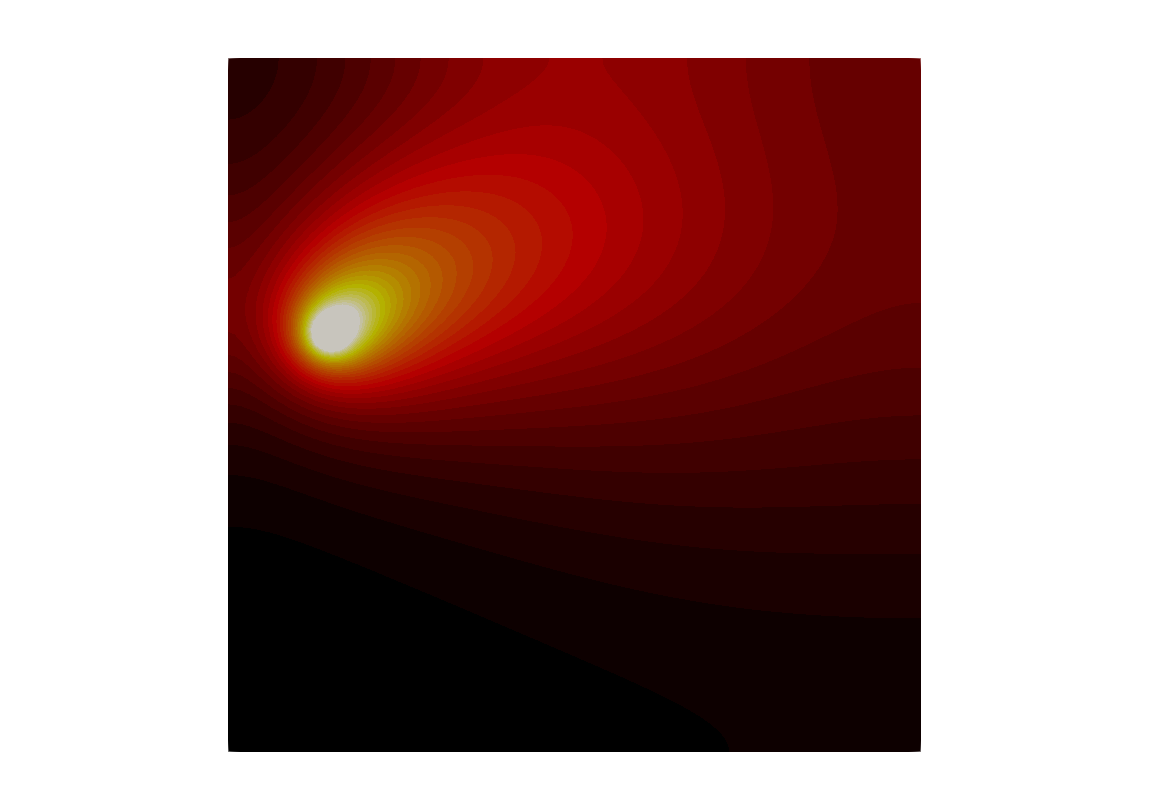}
%}
\caption{Circular arc scan track: Evolution of the temperature distribution for admissible mesh at time steps 1, 10 and 20 (from left to right).\label{fig:TemperatureEvolution}}
\end{figure}
\par
Finally, we want to investigate the behavior of the proposed algorithm compared to a non-admissible grid and to the function support balancing algorithm with function support balancing parameter set to 1 (for further details see~\citet[Sec.~7]{Lorenzo_2017}). 
To this end, we decrease the refinement D\"orfler parameter $\alpha_r$ to 0.059 in the non-admissible case such that we can obtain a mesh with 8 levels of refinement at each time step also for this grid, which otherwise (keeping $\alpha_r=0.1$) would return a coarser grid than the admissible one using the same number of time steps. 
Obviously, the same result for non-admissible meshes can also be obtained increasing the maximum number of iterations, but we prefer to modify the D\"orfler parameter in an effort to obtain a fair comparison between the different discretizations.

\begin{figure}[t!]
\centering
%\subfloat[\label{Mesh1TSNoBal}]{% 
\includegraphics[ width=0.45\textwidth, trim= 80mm 0mm 80mm 20mm,clip=true]{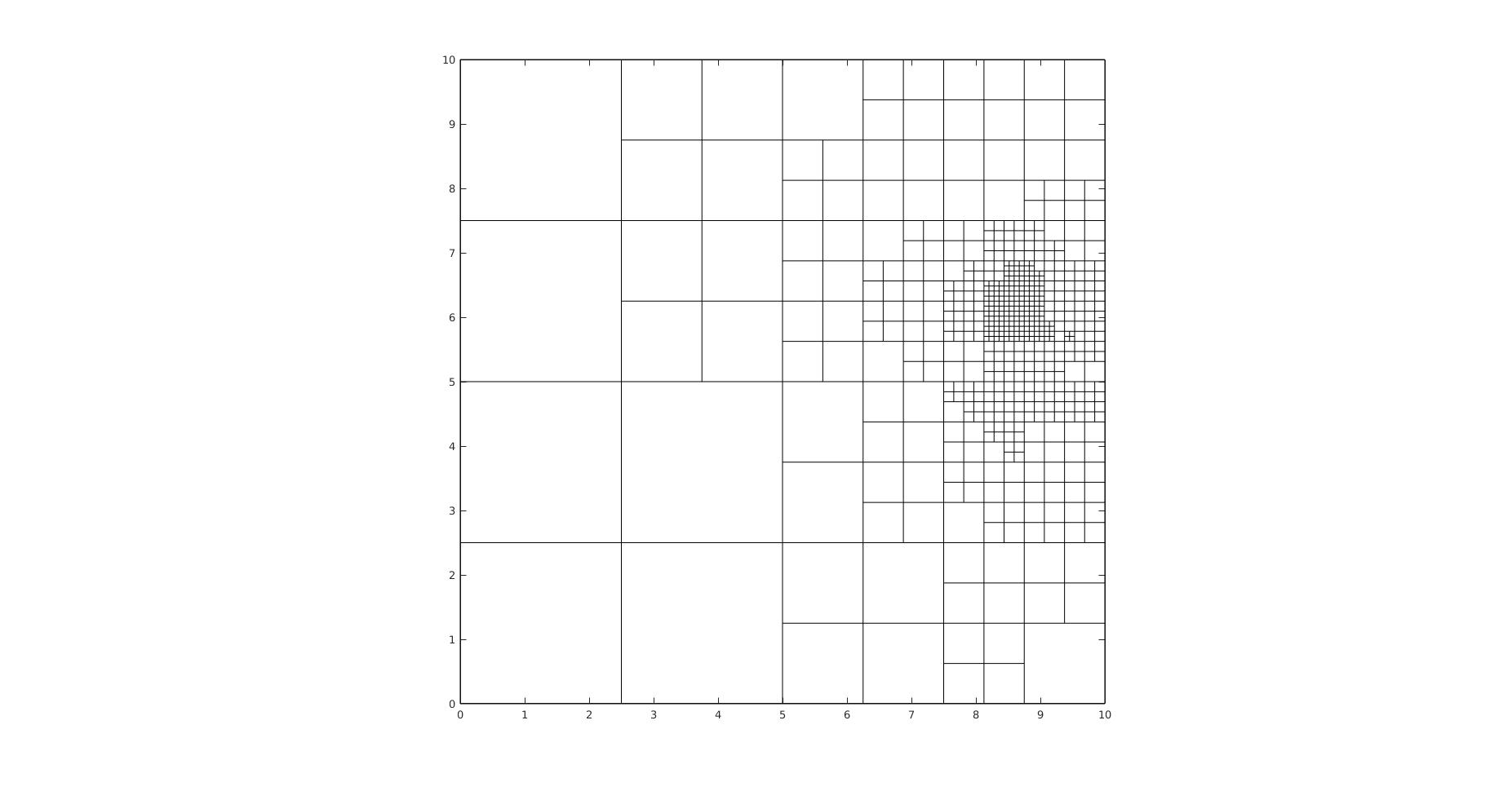}
  %}
%\hfill
%\subfloat[\label{Mesh1TS}]{% 
\includegraphics[ width=0.45\textwidth, trim= 80mm 0mm 80mm 20mm,clip=true]{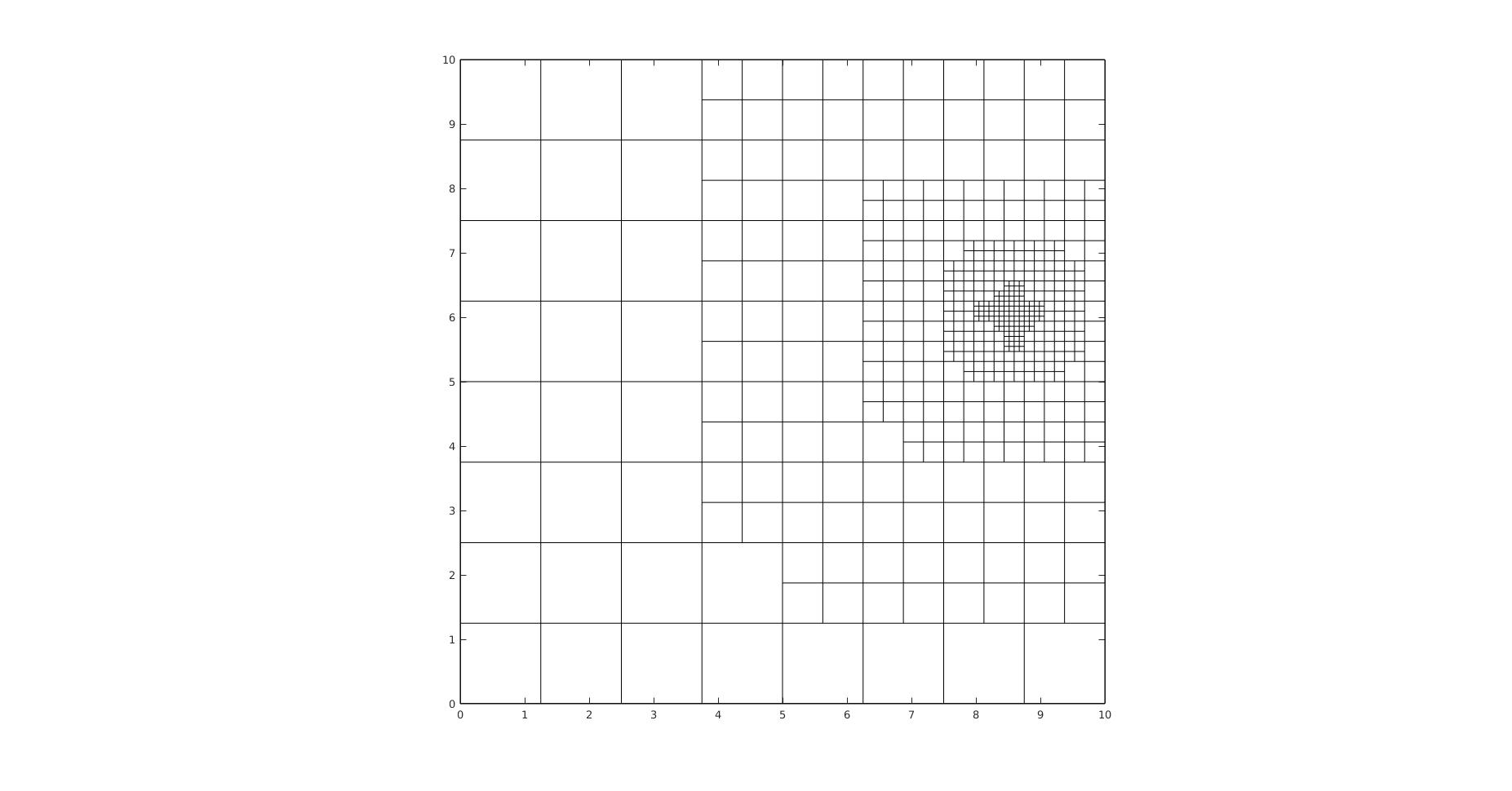}
  %}
%
\linebreak
%\subfloat[\label{Mesh10TSNoBal}]{% 
\includegraphics[width=0.45\textwidth, trim= 80mm 0mm 80mm 20mm,clip=true]{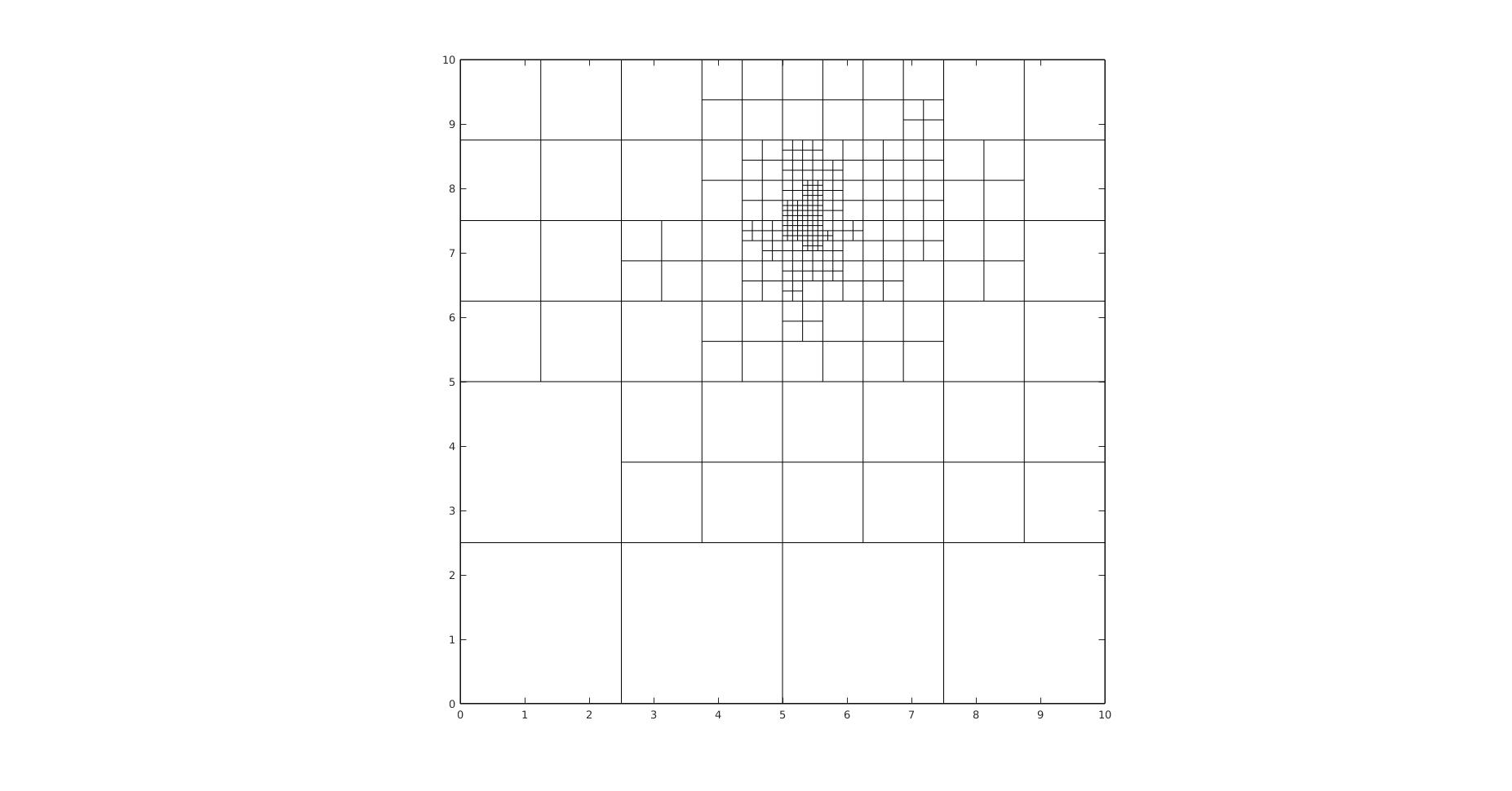}
  %} 
%  \hfill
%\subfloat[\label{Mesh10TS}]{% 
\includegraphics[width=0.45\textwidth, trim= 80mm 0mm 80mm 20mm,clip=true]{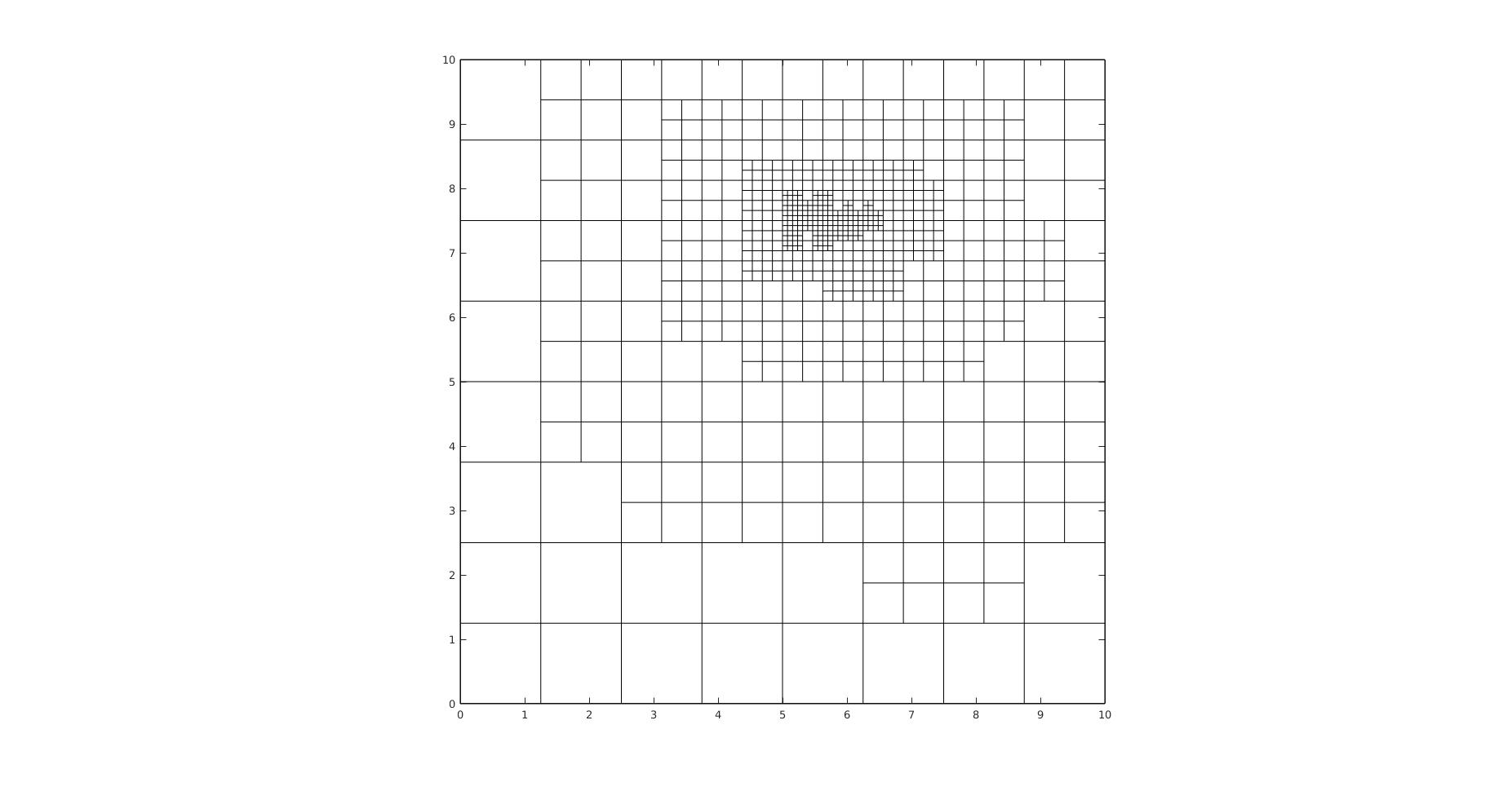}
  %} 
%
\linebreak
%\subfloat[\label{Mesh20TSNoBal}]{% 
\includegraphics[width=0.45\textwidth, trim= 80mm 0mm 80mm 20mm,clip=true]{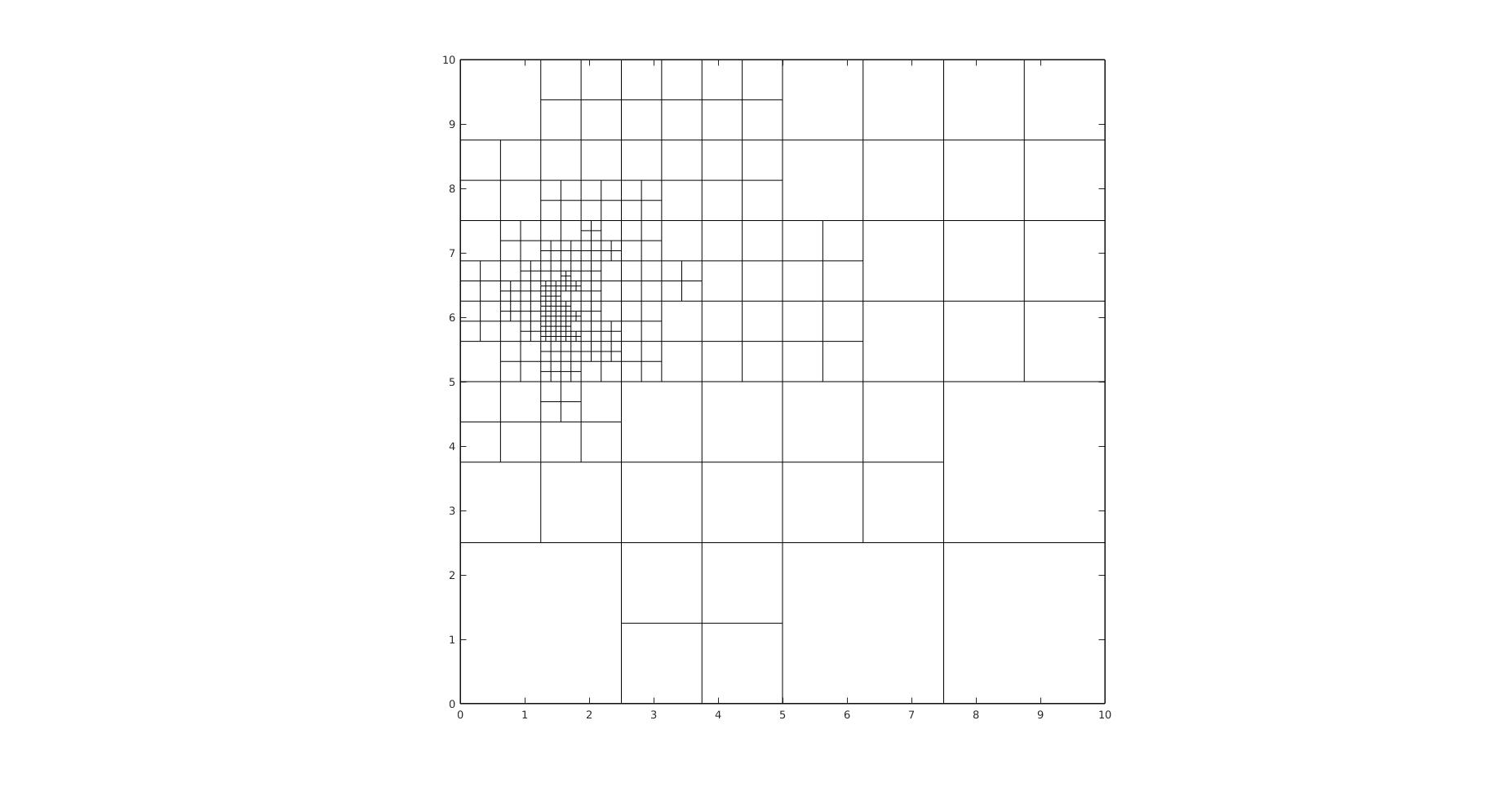}
%}
%\hfill 
%\subfloat[\label{Mesh20TS}]{% 
\includegraphics[width=0.45\textwidth, trim= 80mm 0mm 80mm 20mm,clip=true]{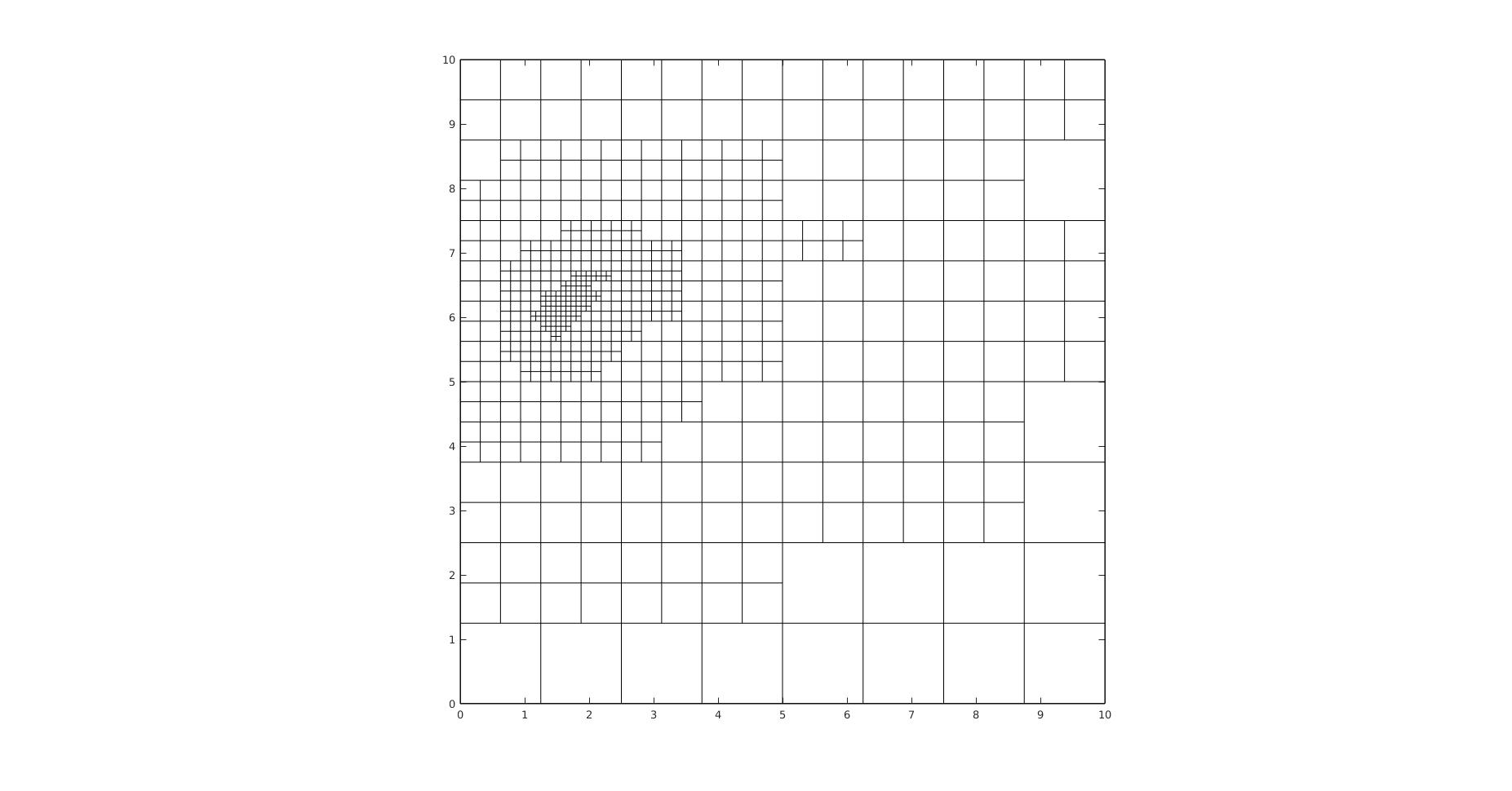}
%}
\caption{Circular arc scan track: Evolution of the non-admissible (left) and admissible (right) adaptive meshes at time steps 1, 10, 20 (from top to bottom).\label{fig:MeshEvolution}}
\end{figure}

\cref{NotAdmComparisonPlots} presents the comparison between the three different cases in terms of total number of DOFs and relative error at each time step, respectively.
We can observe that both the function support balancing and the admissible adaptive discretization reach the same level of accuracy, while the non-admissible grid has a non constant behavior and presents in many time steps a much higher error. 
\cref{fig:TemperatureEvolution} reports the temperature distribution at different time steps using admissible grids, in this case we cannot graphically observe any substantial difference between the admissible and non-admissible results which is instead captured when we look at the relative energy error of~\cref{NotAdmComparisonPlots}.
In any case we can observe how the steep temperature gradients in the proximity of the laser spot can be optimally captured by means of high-order and highly continuous approximation schemes, while the locality of the solution naturally calls for an adaptive discretization.
The comparison of~\cref{fig:MeshEvolution} between admissible and non-admissible grids shows that we obtain much more graded meshes if admissibility requirements are matched.
We can conclude that the proposed algorithm leads to an optimal trade-off between accuracy and number of DOFs (and consequently memory consumption and computational efficiency). 
In fact, we obtain an error comparable with the uniform grid drastically reducing the total number of DOFs, while the function support balancing algorithm leads to a similar accuracy but with almost two times more DOFs per time step.
Finally, we want to remark that we observe a similar behavior in terms of DOFs per time step also for  hierarchical meshes with a higher degree.
%---------------------------------------------------------------------------------------------------------------------------
\subsection{Alternate scan directions} \label{sec:AltScans}
\begin{figure}[h!]
\centering
\includegraphics[width=0.45\textwidth]{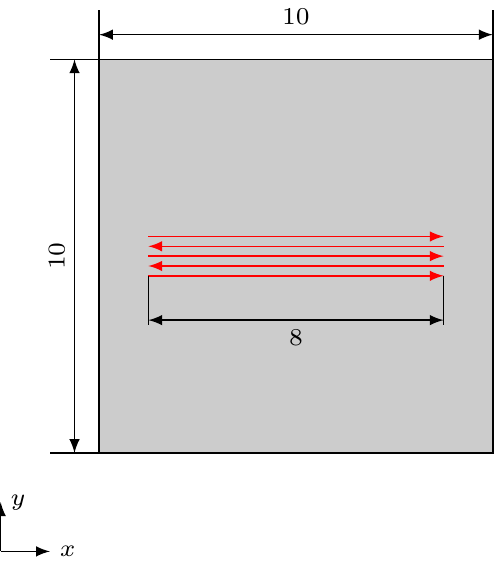}
\caption{Alternate scan directions: Heat source moving on multiple adjacent tracks (in red) with alternate directions.\label{MultiTrackPathSetup}}
\end{figure}

In this second example we consider a moving heat source traveling along multiple, adjacent, 8 mm long tracks in alternate directions on a surface of $10\times 10$ mm$^2$, as depicted in~\cref{MultiTrackPathSetup}.
The external heat source is again defined using the Gaussian distribution of the previous case~\eqref{rhsThermal}, but here we consider a two times smaller radius. 
This scale ratio (between the domain and the heat source radius) is close to the typical one we can find in L-PBF applications. 
We set the laser scan distance between two consecutive tracks (hatch distance) equal to the laser radius, this is a typical choice in L-PBF processes since it avoids gaps between solidified material regions.
For this example we set $\alpha_r=0.08$ and $\alpha_c=0.25$, while the other problem parameters are reported in~\cref{MTParameterTable}.
\begin{table}
\centering
\begin{tabular}{ | c | c | } 
\hline
\textbf{Parameters} & \textbf{Values} \\ 
\hline
Laser power P & $190$[W] \\ 
\hline
Laser speed & $8\left[\text{mm}/\text{sec}\right]$ \\ 
\hline
Absorptivity $\eta$ & $0.33$ \\ 
\hline
Laser radius $r_h$ & $50\left[\mu\text{m}\right]$ \\
\hline
Conductivity $k$ & $29.0\times 10^{-3}\left[\text{W}/\text{m} / \text{K}\right]$ \\
\hline
Specific heat capacity $C_p$ & $650.0\left[\text{J}/\text{kg}/\text{K}\right]$ \\
\hline
Density $\rho$ & $8440.0\left[\text{kg}/\text{m}^3\right]$ \\
\hline
Initial temperature $\theta_0$ & $25.0\left[^{\circ}\text{C}\right]$ \\
\hline
Hatch distance & $50\left[\mu\text{m}\right]$ \\
\hline
\end{tabular}
\caption{Alternate scan directions: Process and material parameters.\label{MTParameterTable}}
\end{table}
\par
\begin{figure}[h!]
\centering
%\subfloat[\label{TimeMT}]{% 
\includegraphics[width=0.45\textwidth]{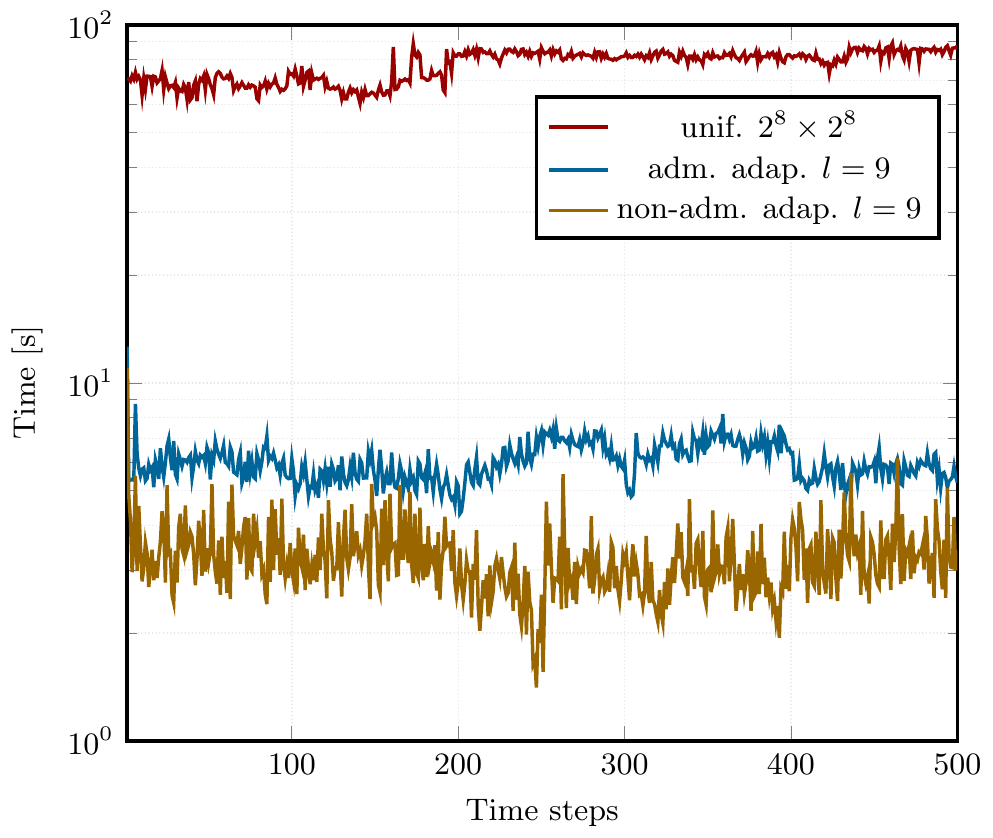}
%  } 
%
\hfill 
%\subfloat[\label{DofsMT}]{% 
\includegraphics[width=0.45\textwidth]{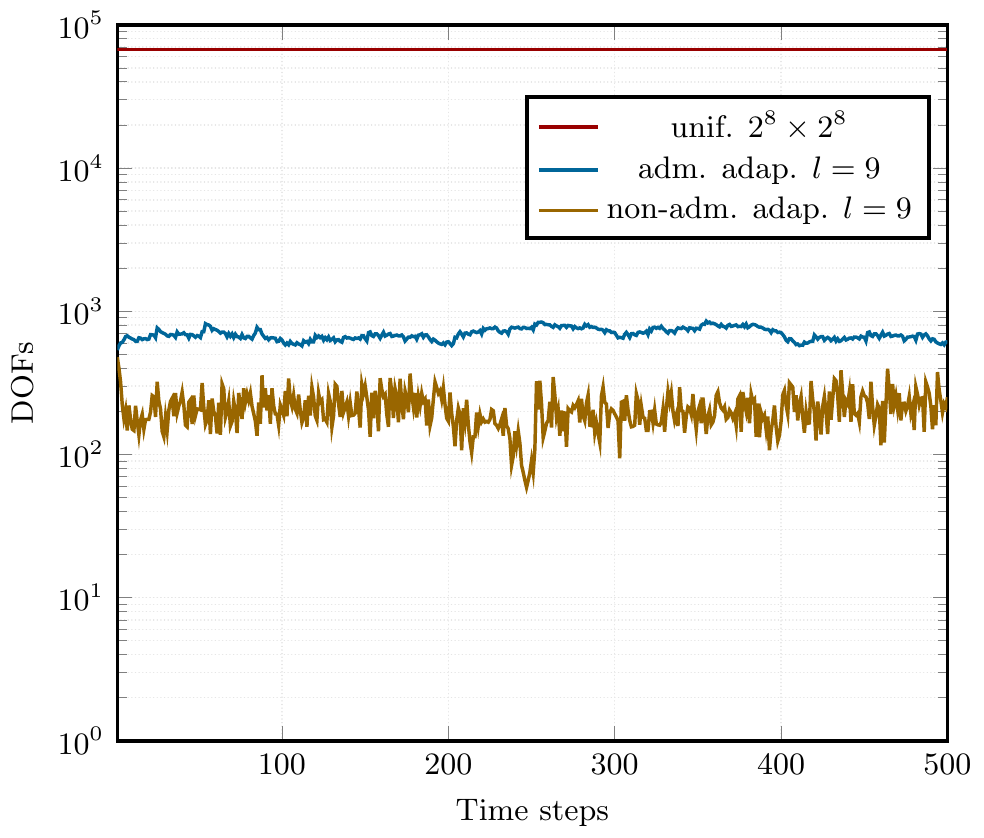}
%}
\caption{Alternate scan directions: CPU time (left) and DOFs  (right)  comparison at each time step between uniform, admissible and non-admissible adaptive meshes with 9 levels of refinements.}
\label{Exm2CPU}
\end{figure}
\begin{figure}[h!]
\centering
%\subfloat[\label{InternalEnergyMT}]{% 
\includegraphics[width=0.45\textwidth]{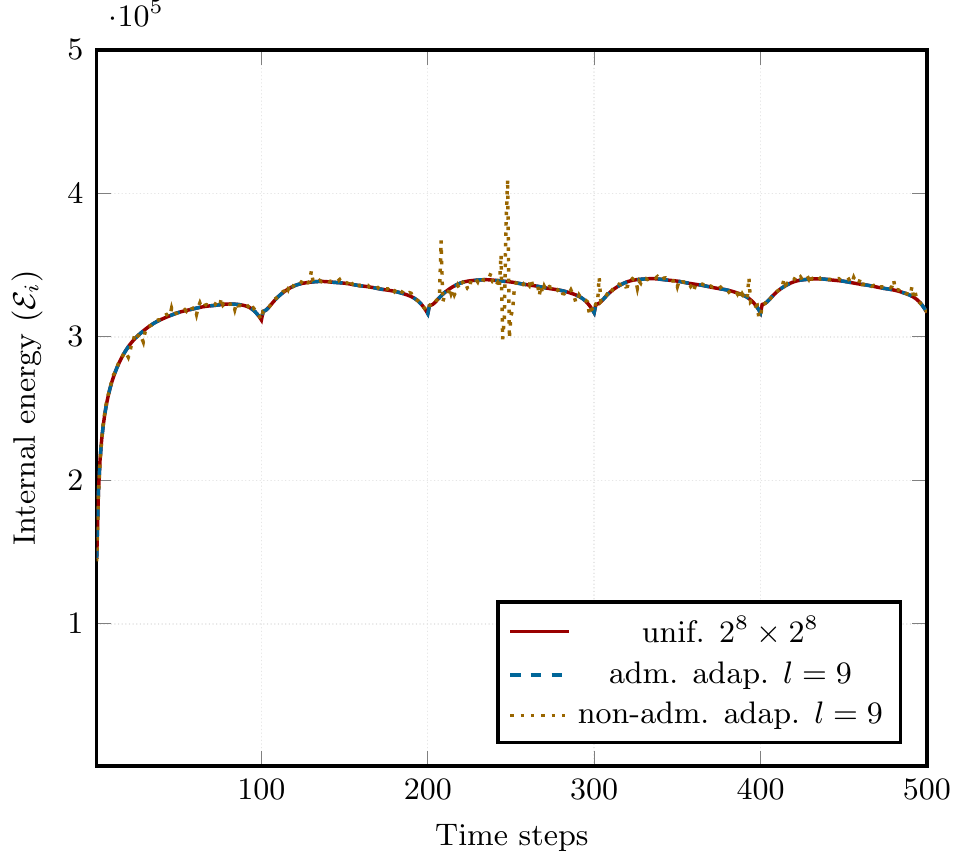}  %} 
\hfill 
%\subfloat[\label{TotalEnergyMT}]{% 
\includegraphics[width=0.45\textwidth]{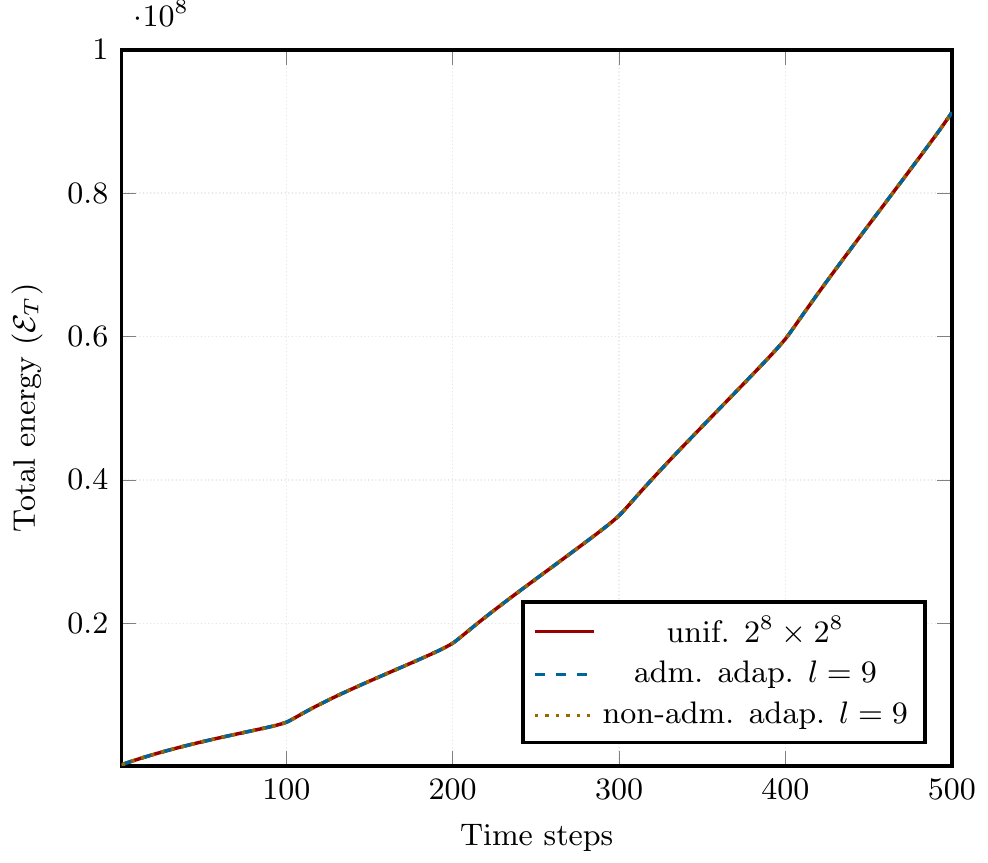}  %} 
\caption{Alternate scan directions: Internal (left) and total (right) energy of the system at each time step between uniform, admissible and non-admissible adaptive meshes with 9 levels of refinements.\label{ErrorMT}}
\end{figure}

We now compare an admissible adaptive mesh with 9 levels of refinement (adm. adap. $l=9$) with respect to a non-admissible adaptive mesh with the same maximum level of refinement (non-adm. adap. $l=9$) and the corresponding uniform mesh (unif. $2^8\times 2^8$).
Analogously to what we did in the previous section, we set a different D\"orfler parameter for the non-admissible discretization ($\alpha_r=0.07$).
~\cref{Exm2CPU} shows the CPU time and the total number of DOFs at each time step for three considered grids. 
We can observe that, compared to the uniform mesh, in both cases we obtain a remarkable advantage when employing an adaptive mesh.
This improvement in terms of both DOFs and CPU time is not affecting the quality of the solution.
In fact, as shown in~\cref{ErrorMT}, the values of both the internal and the total energy at each time step do not present any substantial difference between the uniform and the admissible adaptive case.
On the other side, when we adopt a non-admissible mesh, the internal energy value strongly oscillates from the reference value, while we have a limited advantage in terms of both CPU time and number of DOFs compared to the admissible mesh.
\par

\begin{figure}[h!]
\centering
\includegraphics[width=0.45\textwidth, trim= 100mm 50mm 100mm 0mm,clip=true]{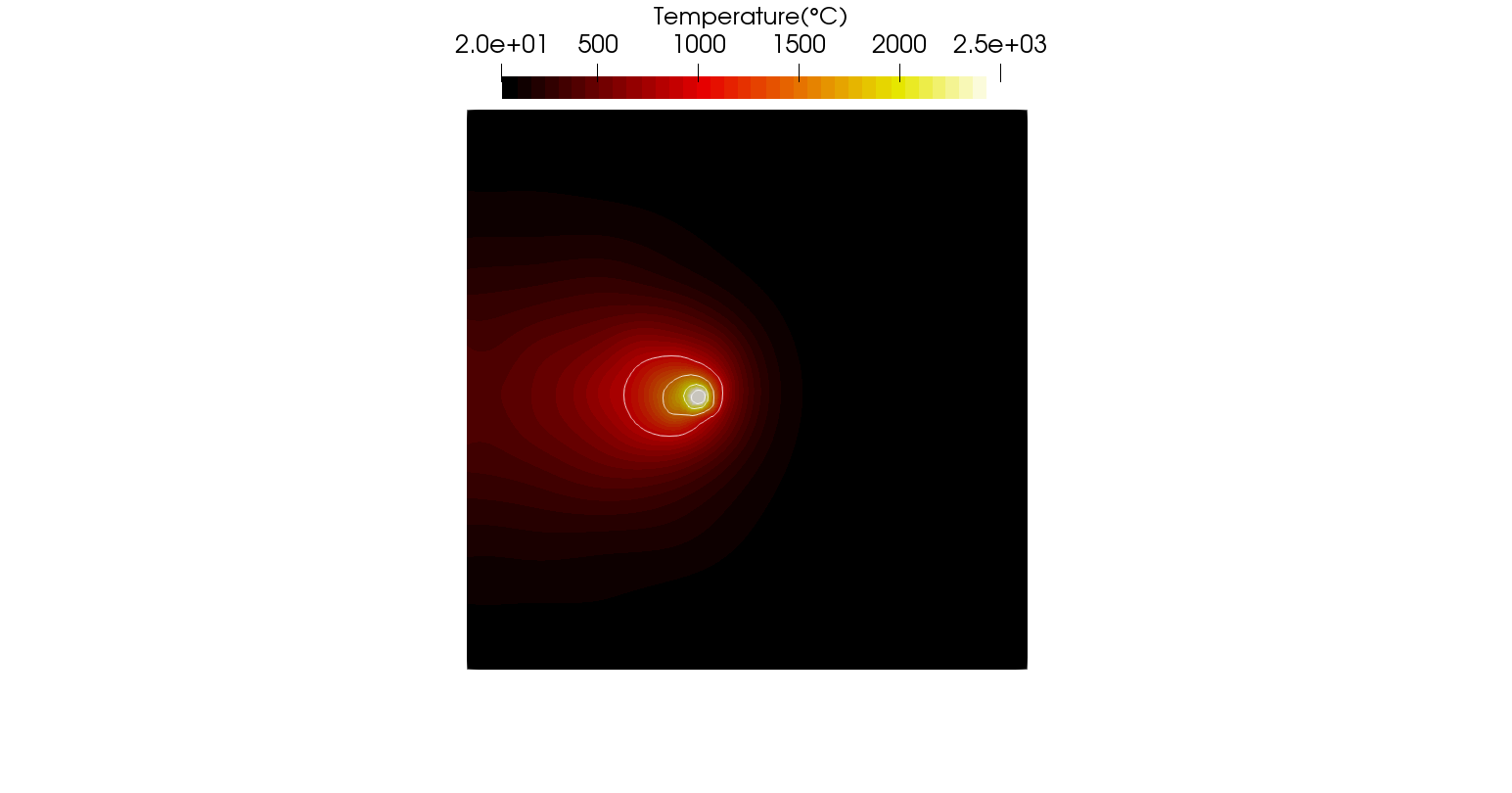}
%}
%
%\hfill
%\subfloat[]{% 
\includegraphics[width=0.45\textwidth, trim= 100mm 50mm 100mm 0mm,clip=true]{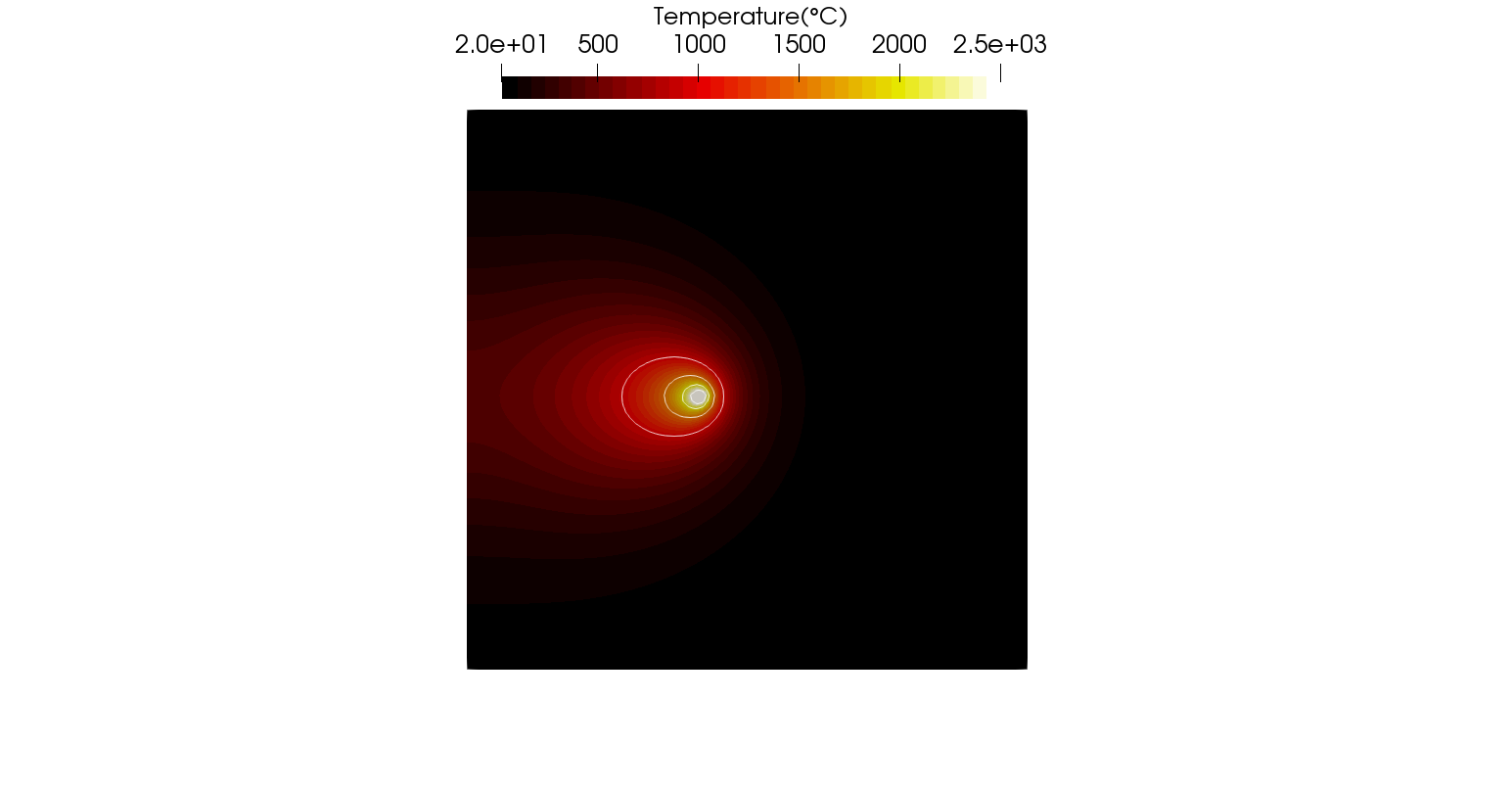}
%}
%
\linebreak
%\subfloat[]{% 
\includegraphics[width=0.45\textwidth, trim= 100mm 50mm 100mm 0mm,clip=true]{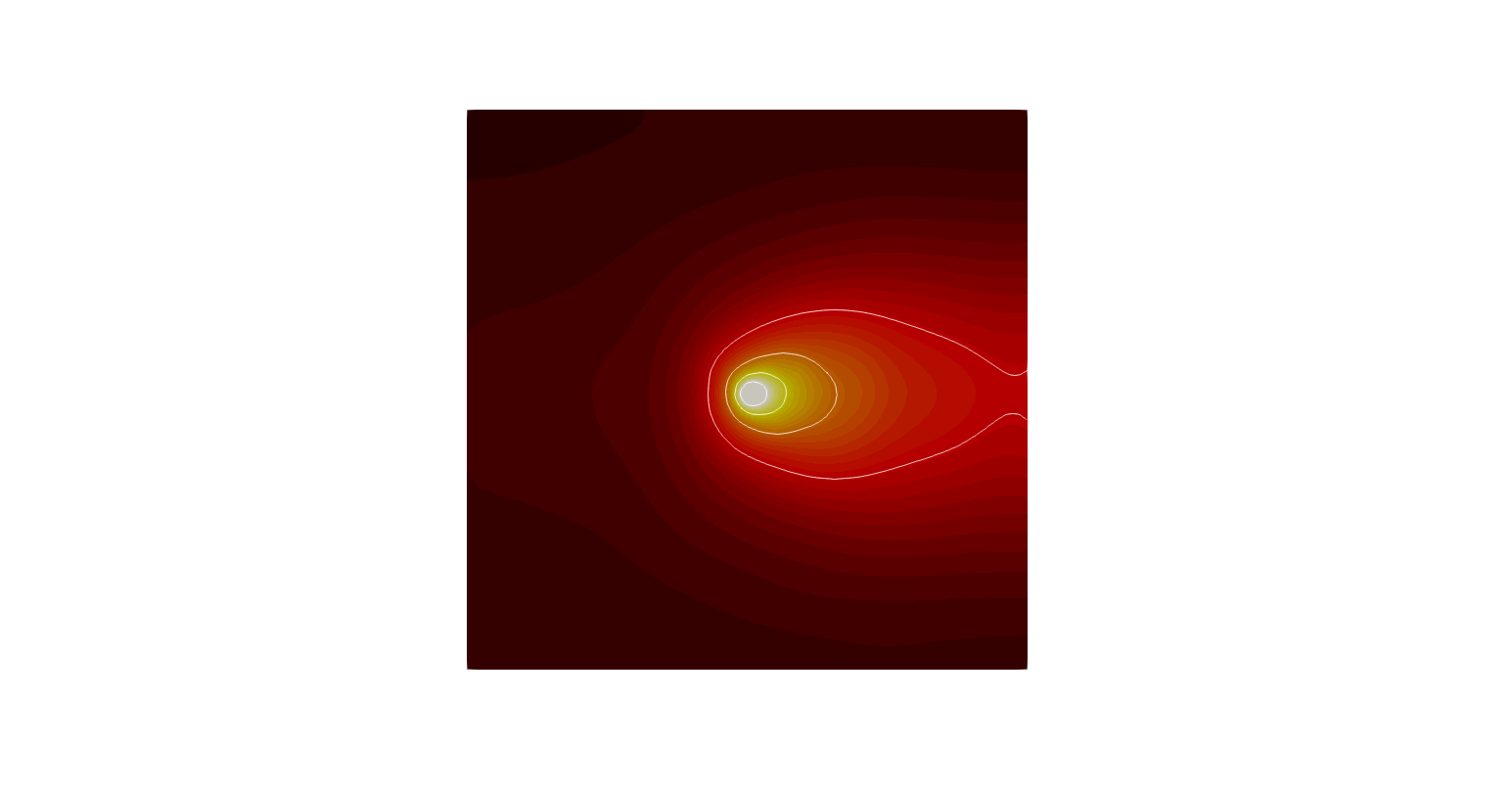}
%}
%
%\hfill 
%\subfloat[]{% 
\includegraphics[width=0.45\textwidth, trim= 100mm 50mm 100mm 0mm,clip=true]{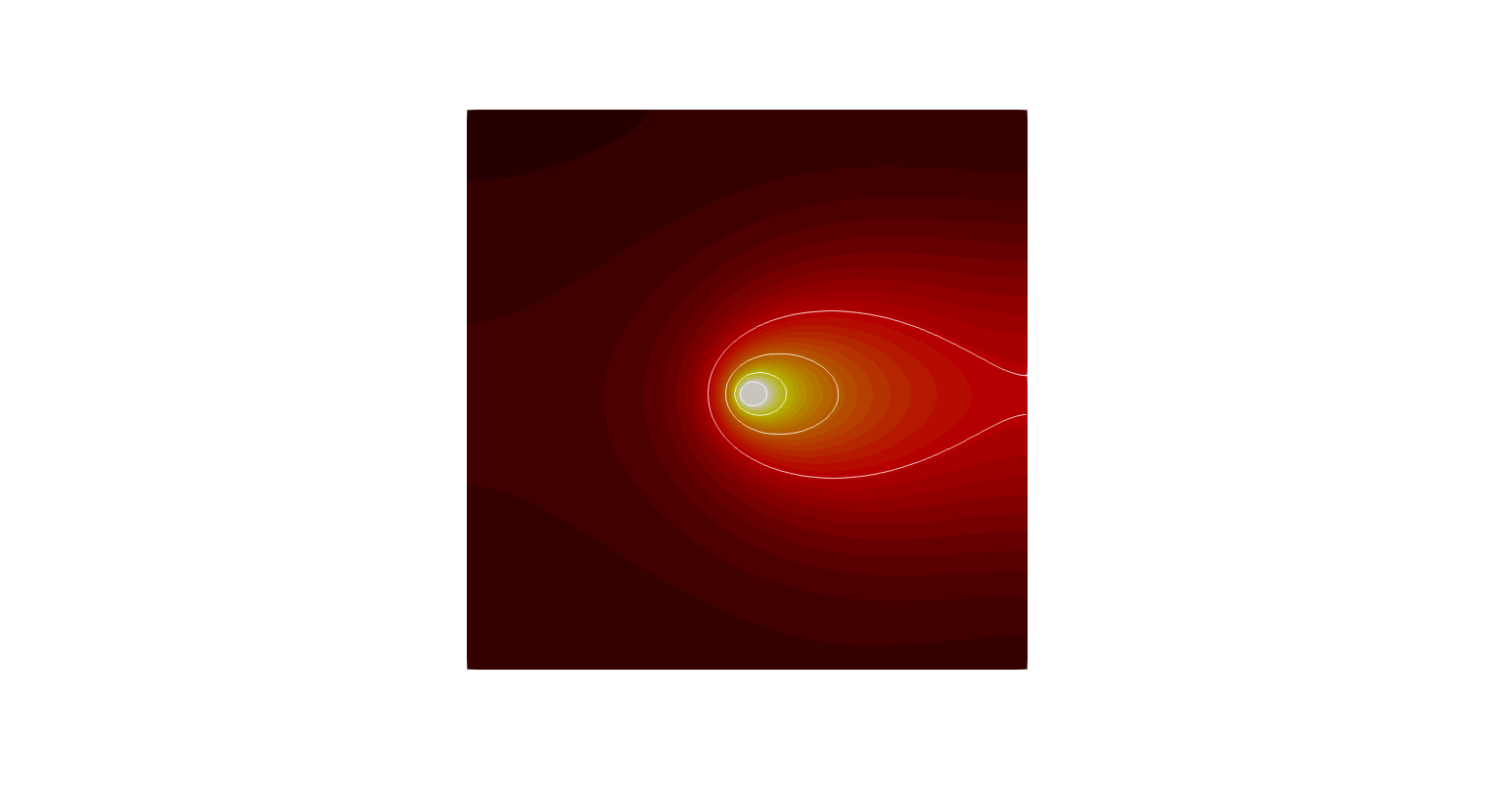}
%}
%
\linebreak
% 
%\subfloat[]{% 
\includegraphics[width=0.45\textwidth, trim= 100mm 50mm 100mm 0mm,clip=true]{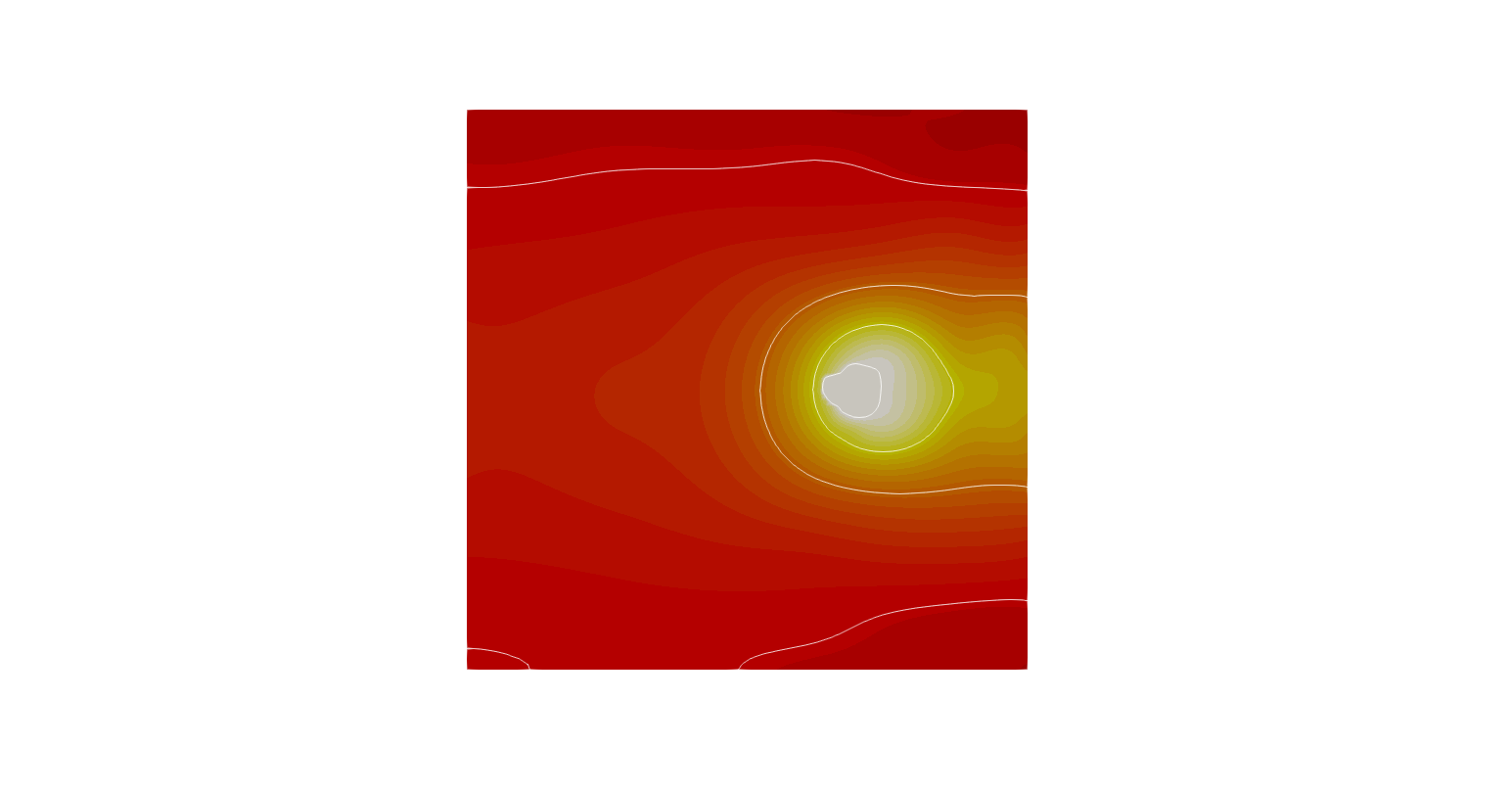}
%} 
%\hfill
%\subfloat[]{% 
\includegraphics[width=0.45\textwidth, trim= 100mm 50mm 100mm 0mm,clip=true]{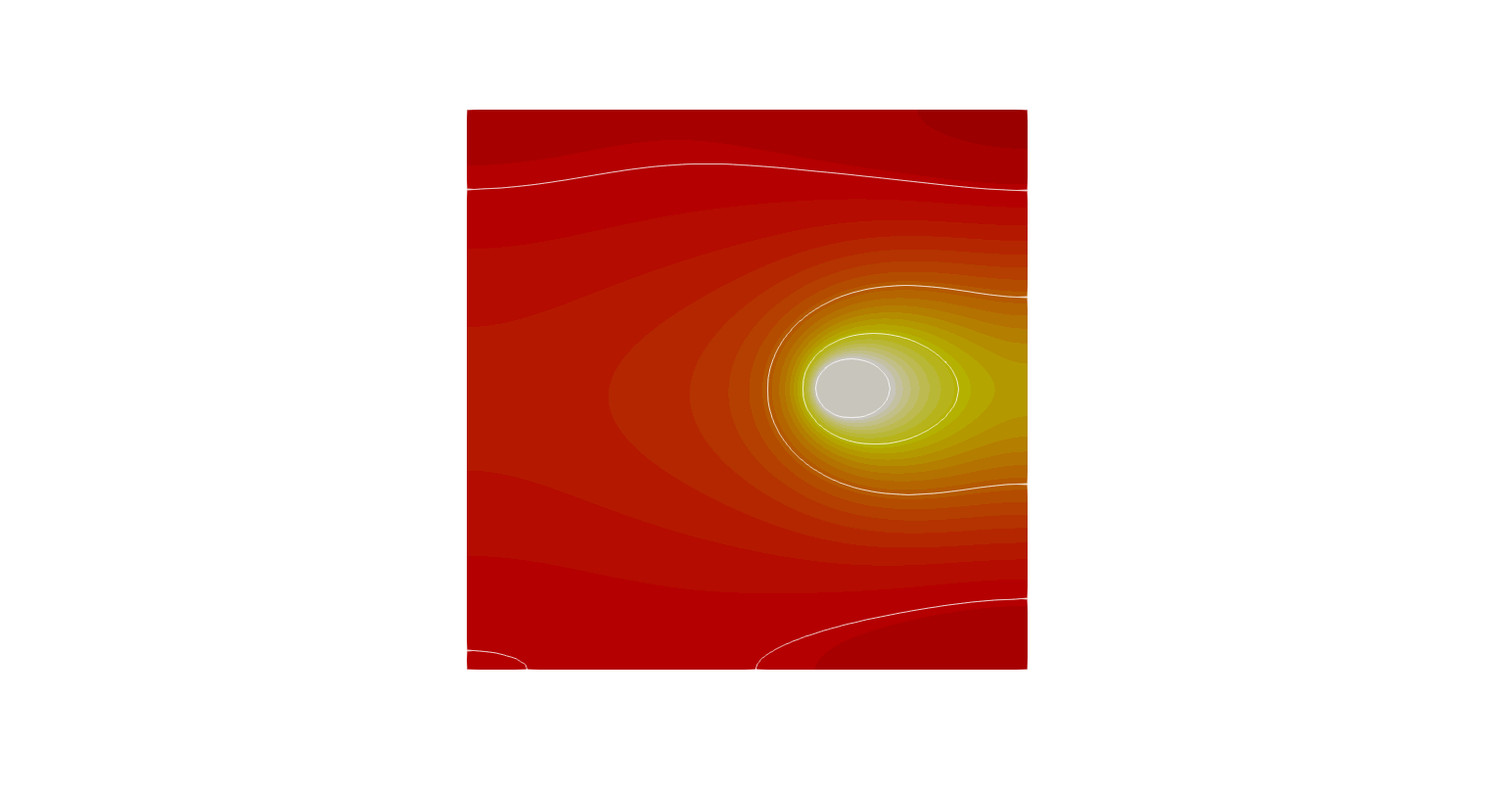}
%} 
%
\linebreak
%
%\subfloat[]{% 
\includegraphics[width=0.45\textwidth, trim= 100mm 50mm 100mm 0mm,clip=true]{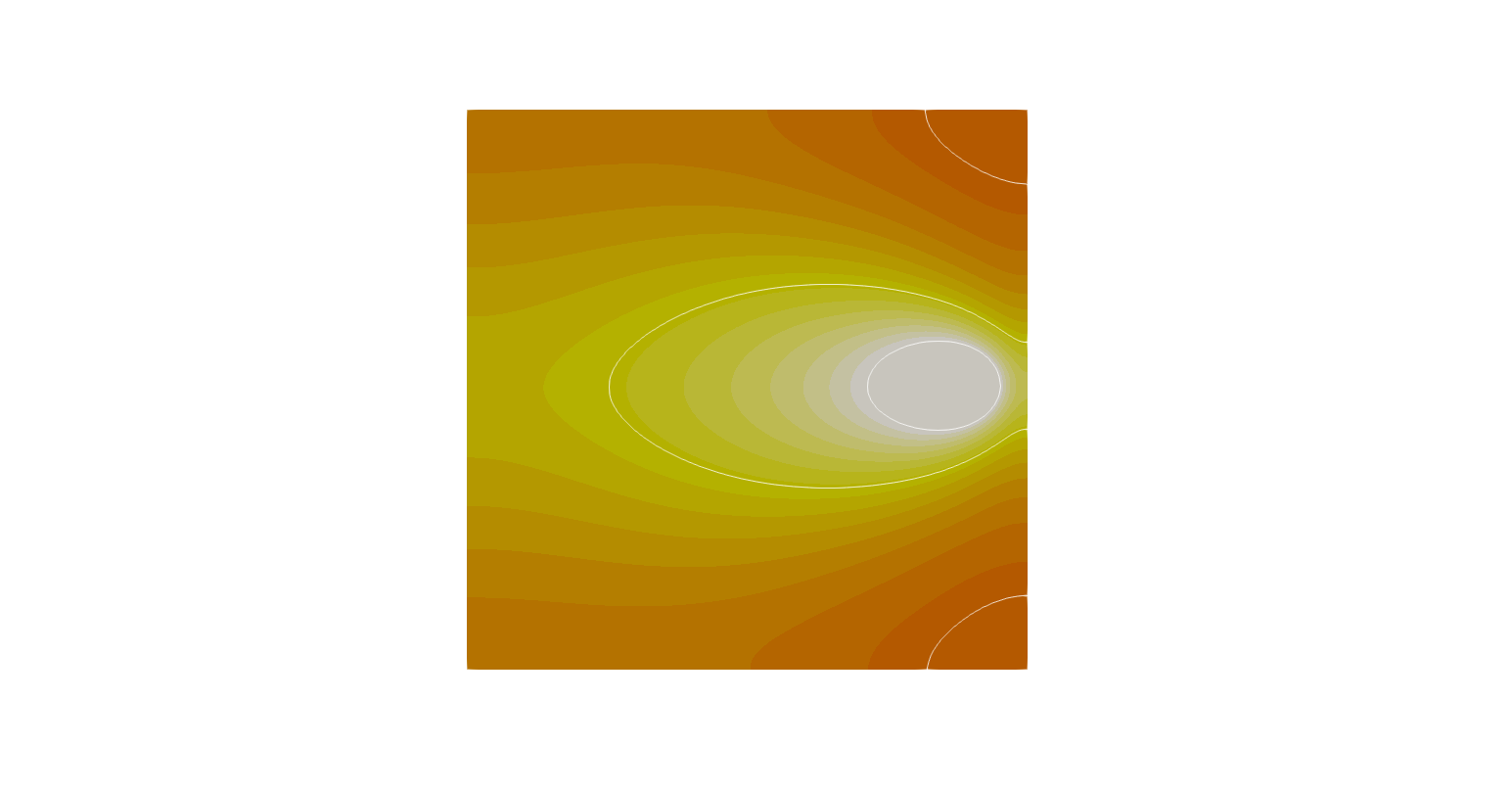}
%}
%\hfill
%\subfloat[]{% 
\includegraphics[width=0.45\textwidth, trim= 100mm 50mm 100mm 0mm,clip=true]{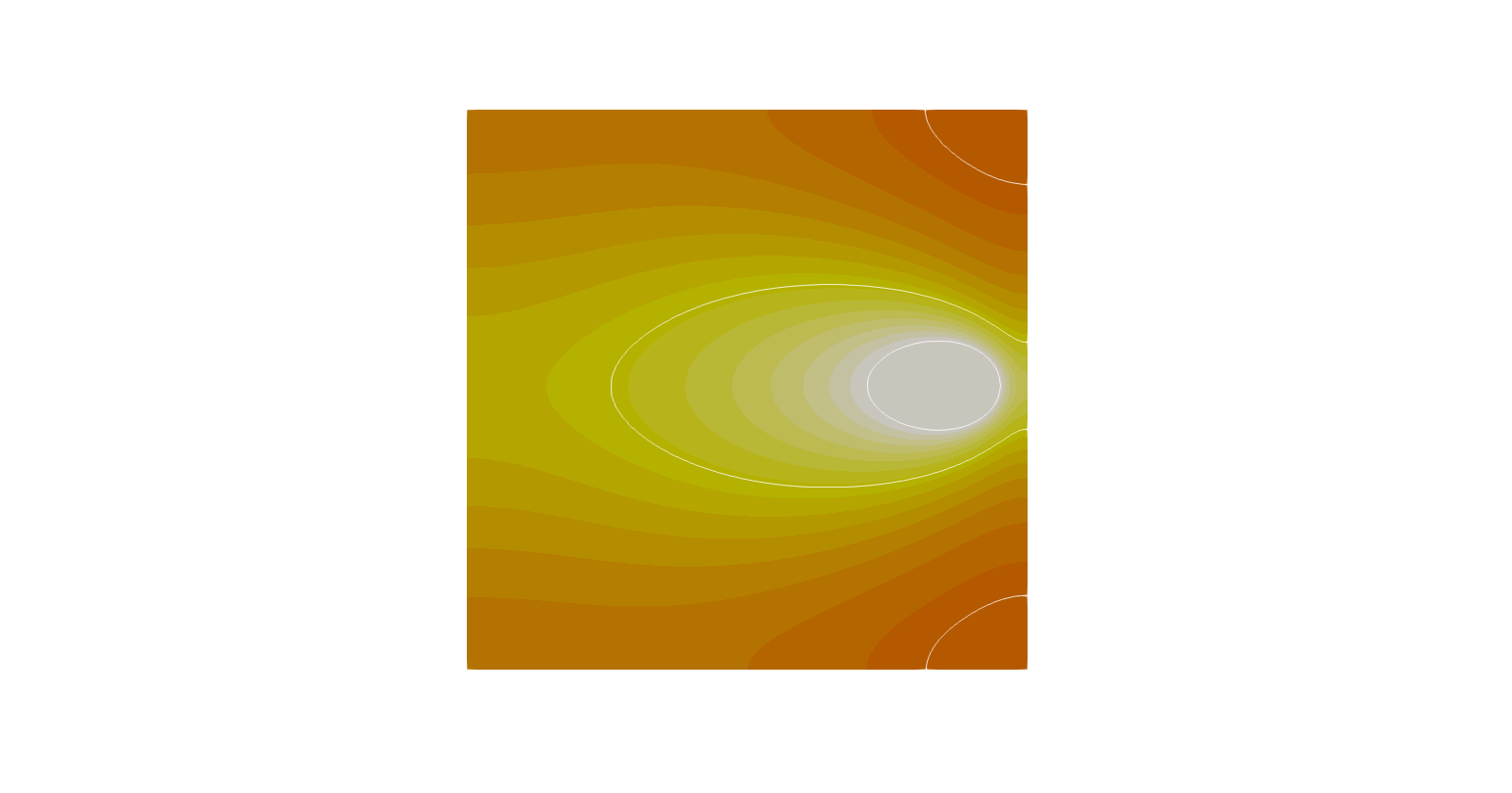}
%}
%
\linebreak
\caption{Alternate scan directions: Evolution of the non-admissible (left) and admissible (right) temperature distributions at time steps 40, 150, 330, and 500 (from top to bottom) with contour lines at 1000, 1500, 2000 and 2500$^{\circ}$C. The complete videos of these simulations are provided together with the supplementary material of this article. \label{fig:MultiTrackTemperatureEvolution}}
\end{figure}

\begin{figure}[h!]
\centering
\includegraphics[width=0.45\textwidth, trim= 80mm 0mm 80mm 20mm,clip=true]{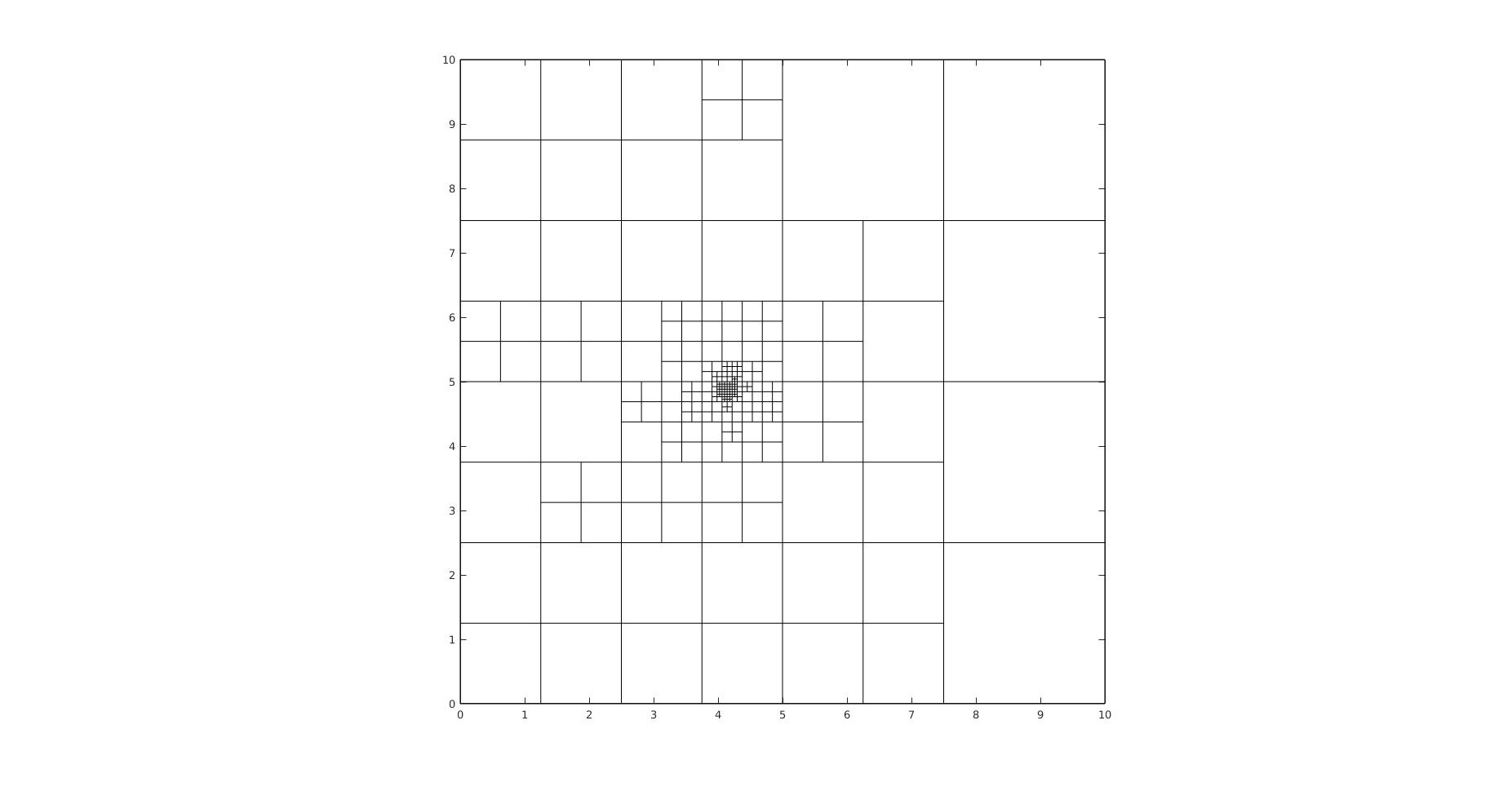}
%\hfill
\includegraphics[width=0.45\textwidth, trim= 80mm 0mm 80mm 20mm,clip=true]{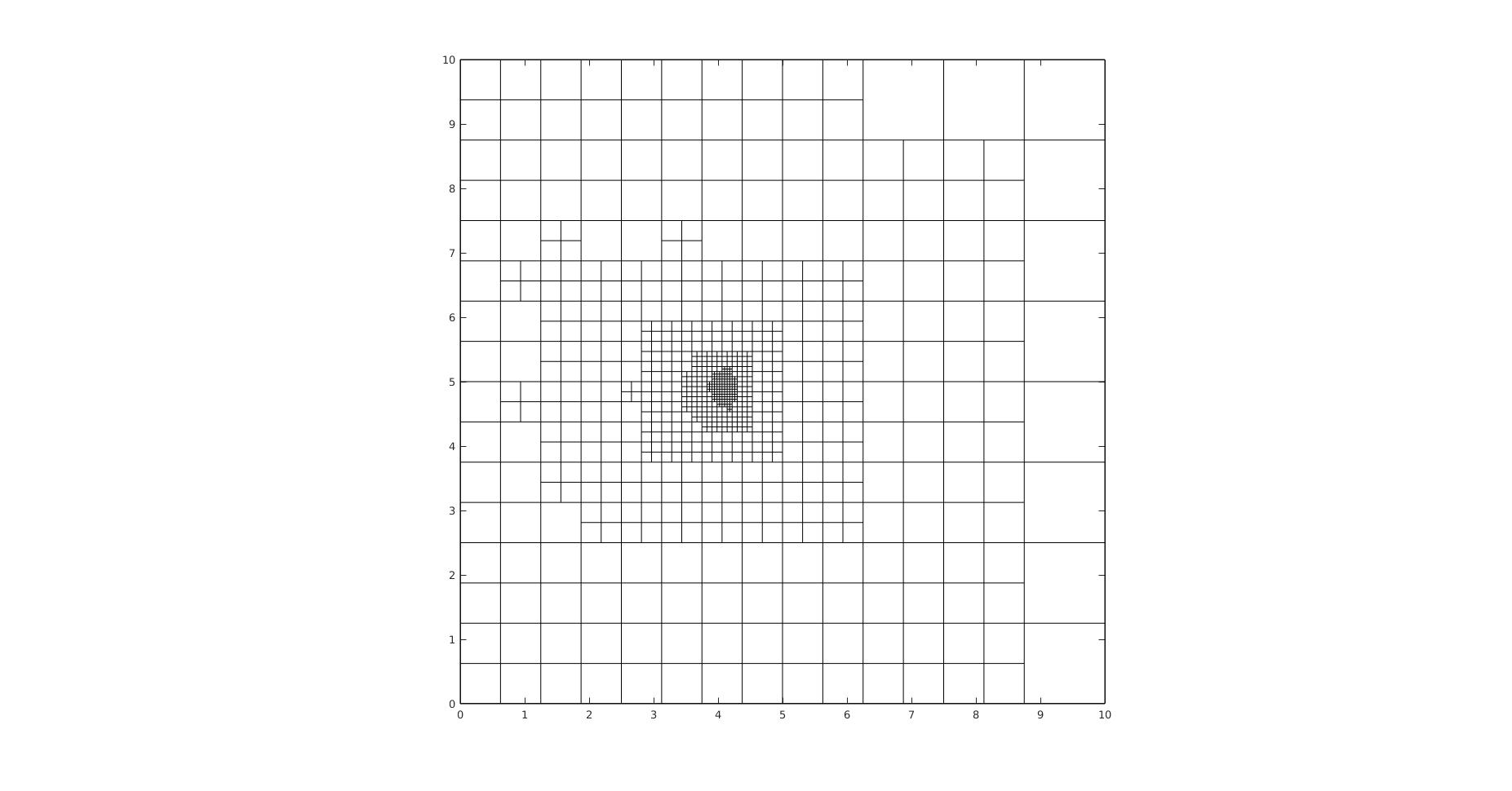}
\linebreak
\includegraphics[width=0.45\textwidth, trim= 80mm 0mm 80mm 20mm,clip=true]{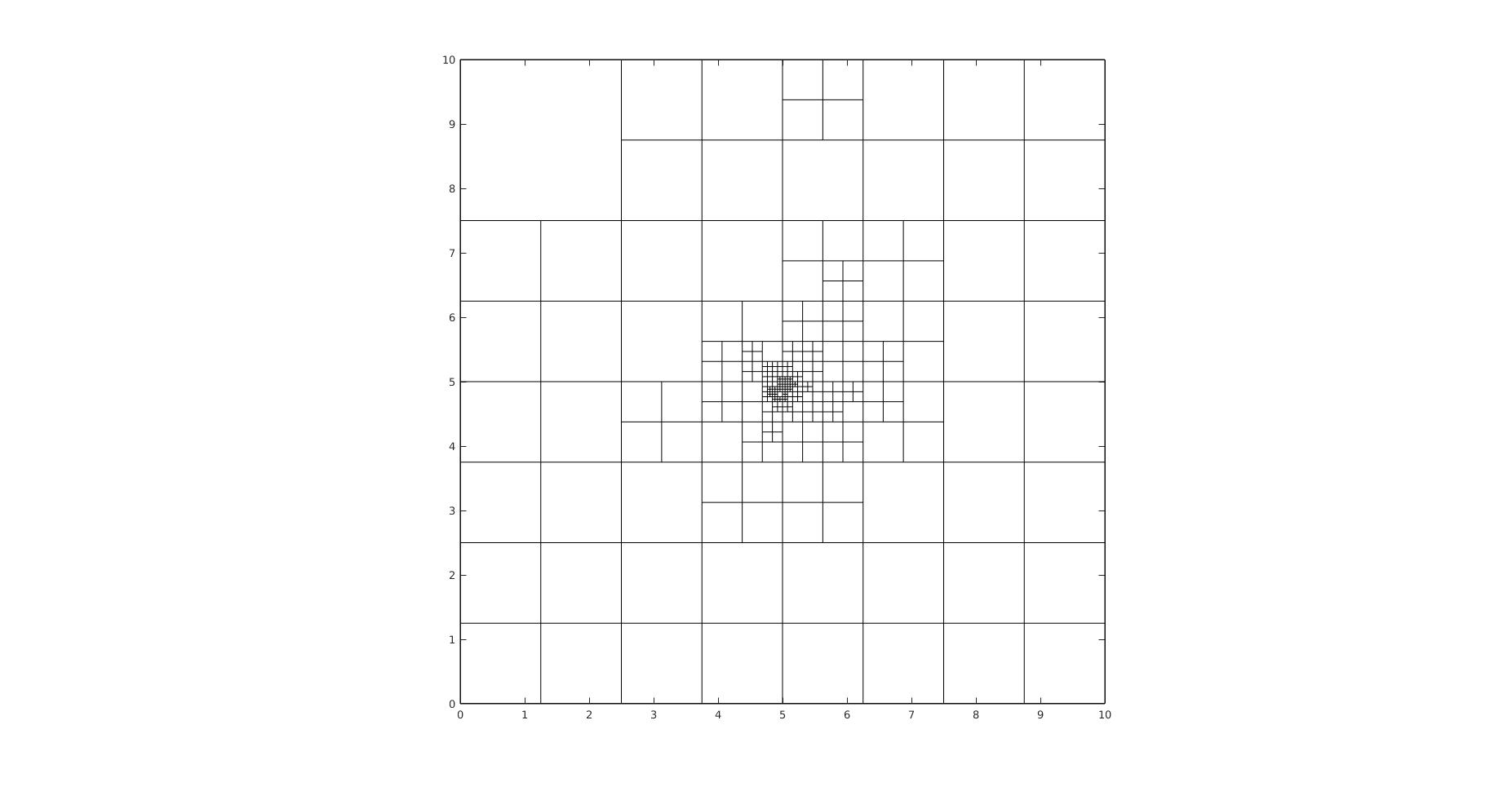}
%\hfill
\includegraphics[width=0.45\textwidth, trim= 80mm 0mm 80mm 20mm,clip=true]{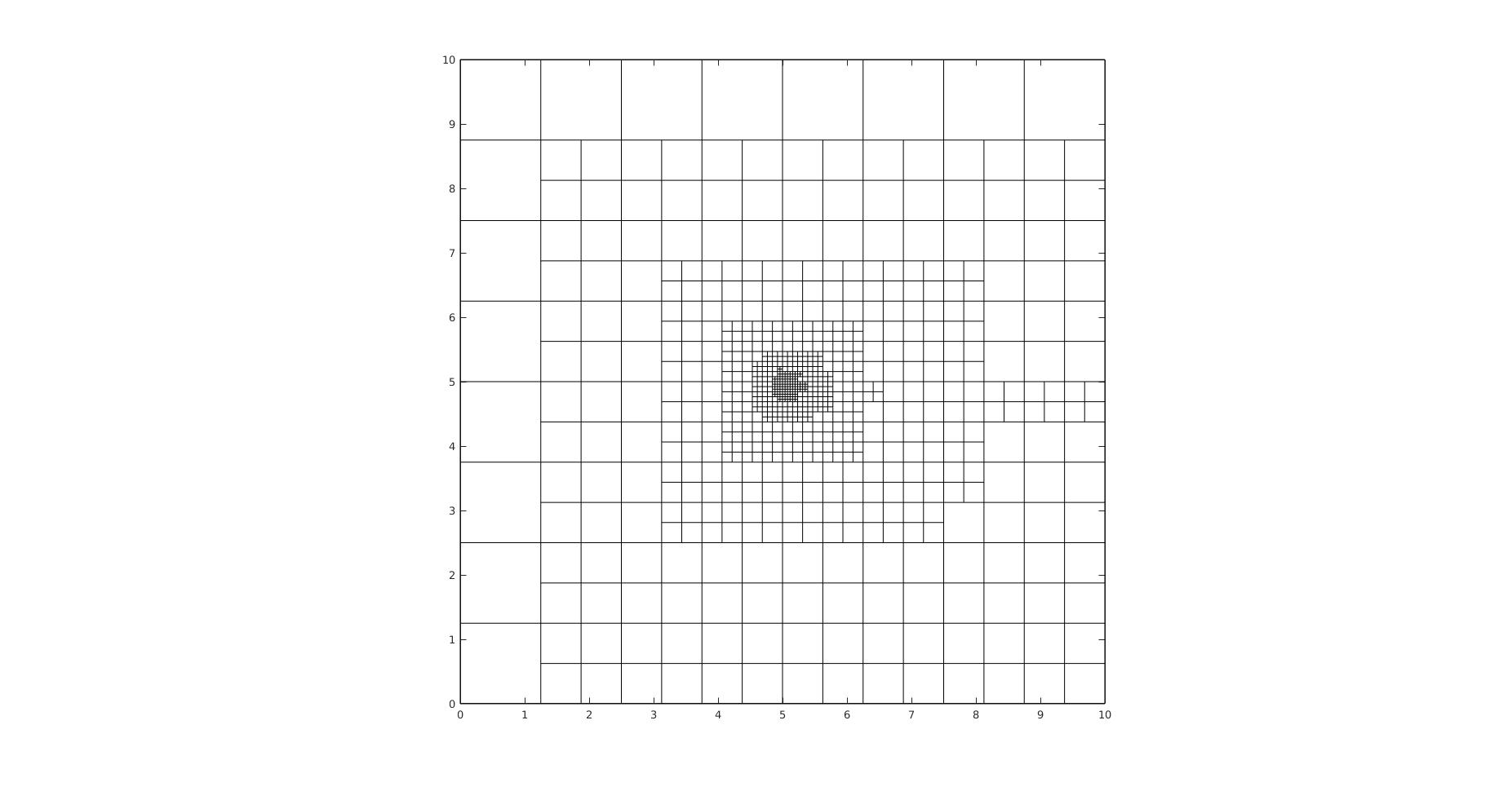}
\linebreak
\includegraphics[width=0.45\textwidth, trim= 80mm 0mm 80mm 20mm,clip=true]{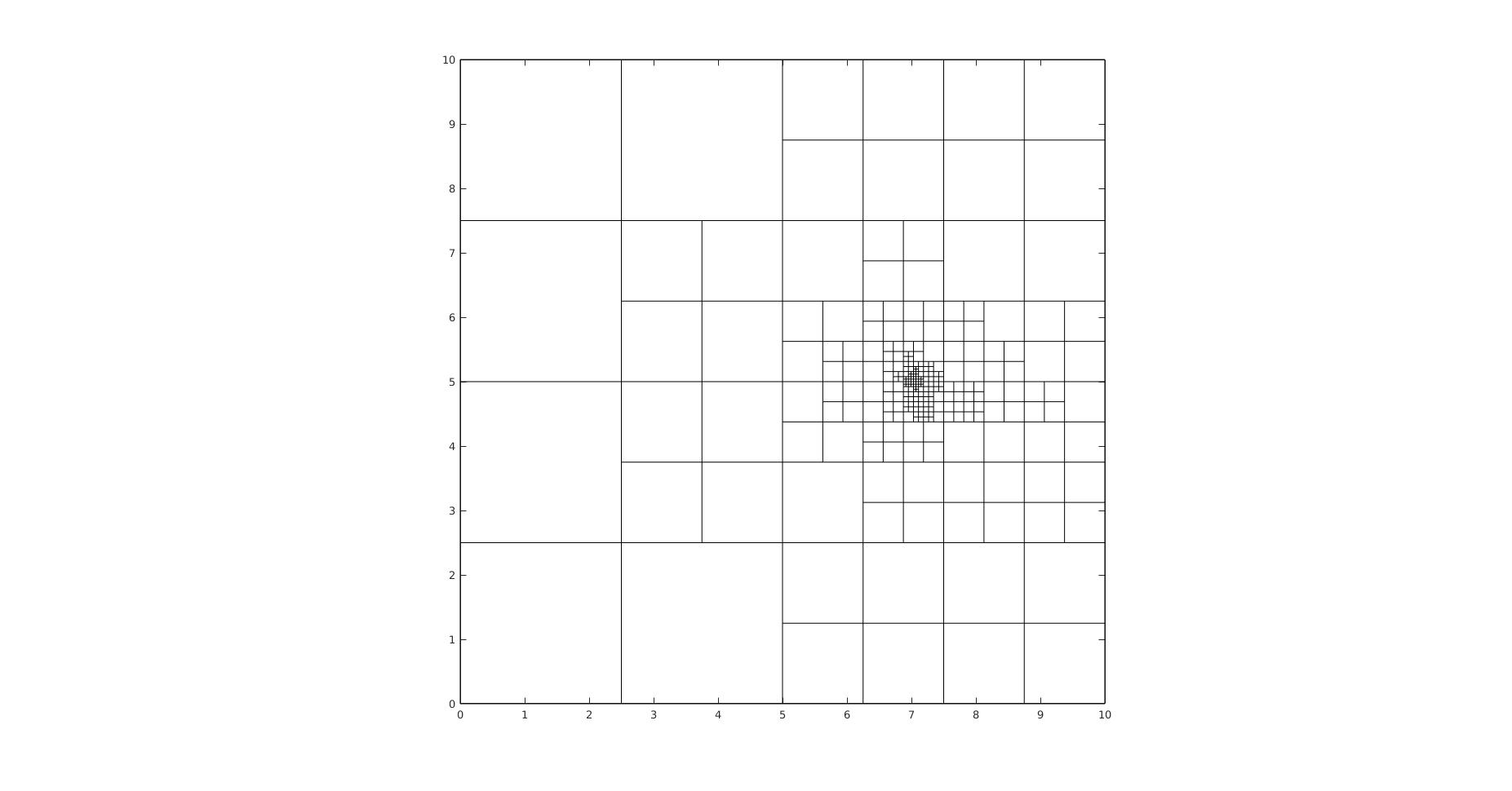}
%\hfill
\includegraphics[width=0.45\textwidth, trim= 80mm 0mm 80mm 20mm,clip=true]{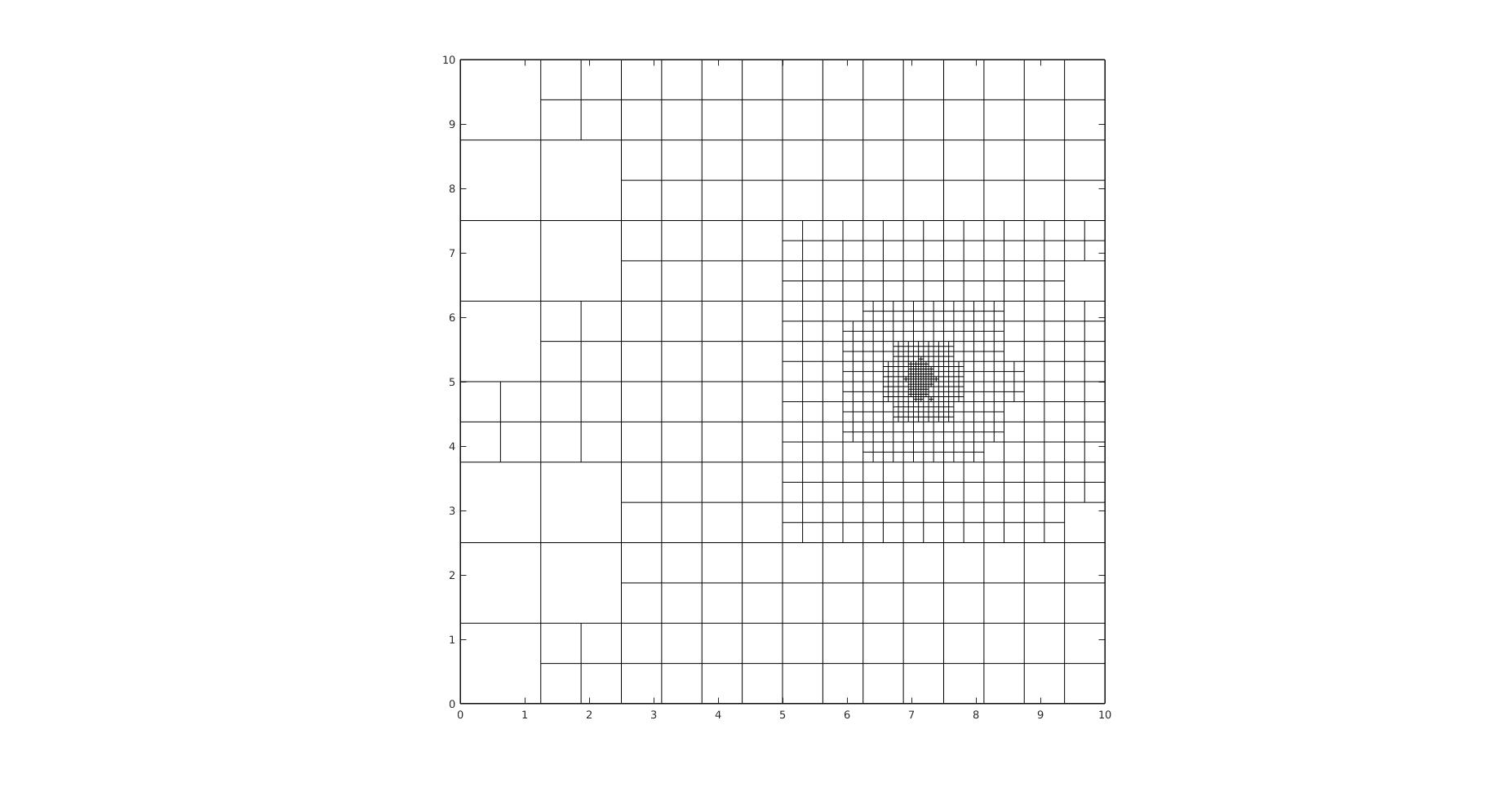}
\linebreak
\includegraphics[width=0.45\textwidth, trim= 80mm 0mm 80mm 20mm,clip=true]{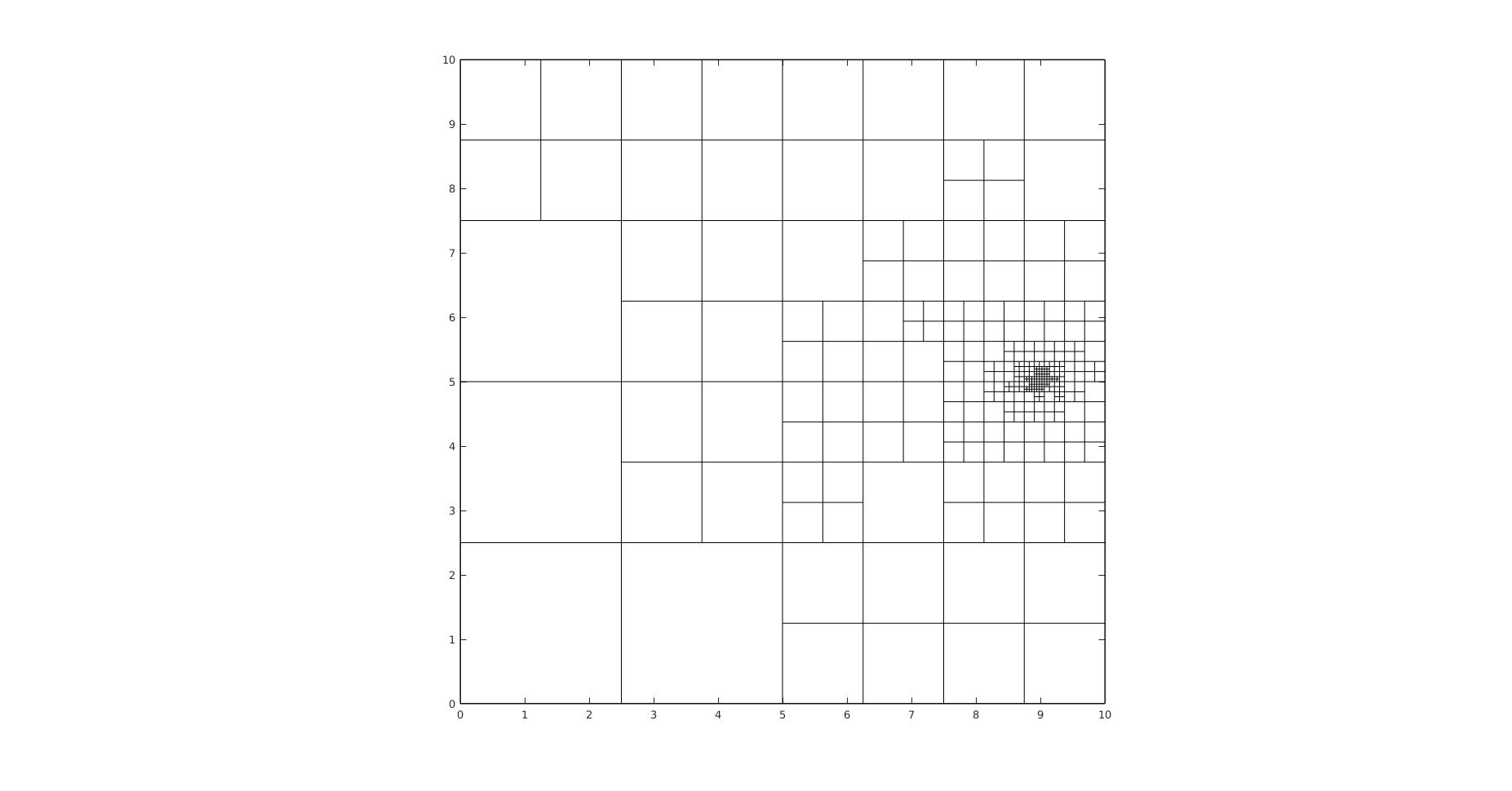}
%\hfill
\includegraphics[width=0.45\textwidth, trim= 80mm 0mm 80mm 20mm,clip=true]{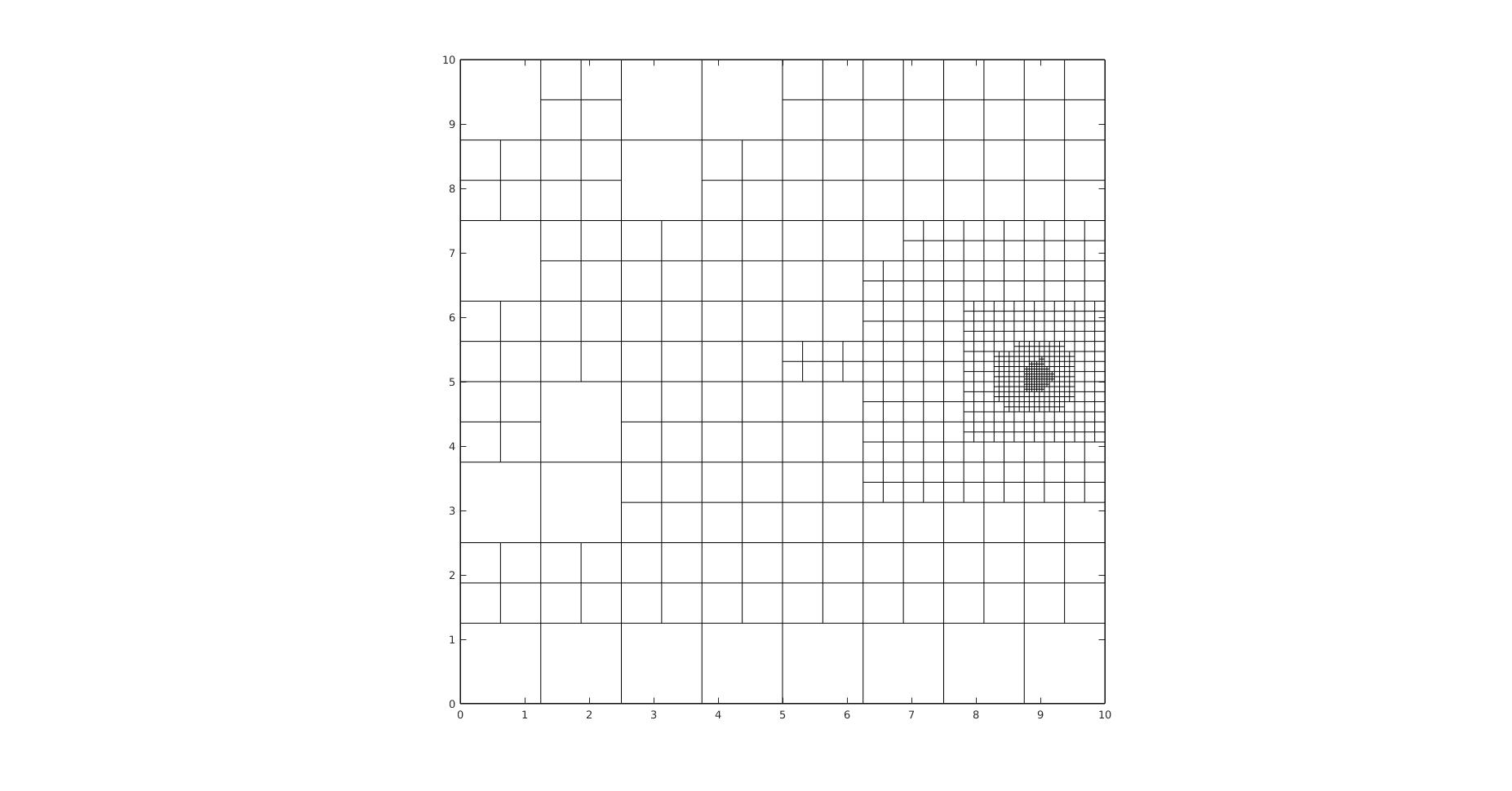}
\caption{Alternate scan directions: Evolution of the non-admissible (left) and admissible (right) adaptive meshes at time steps 40, 150, 330, and 500 (from top to bottom).\label{fig:MultiTrackMeshEvolution}}
\end{figure}

Finally,~\cref{fig:MultiTrackTemperatureEvolution,fig:MultiTrackMeshEvolution} report the  temperature distributions and the corresponding discretizations at different positions of the heat source along the path, comparing admissible and non-admissible results. 
We note that, even if the solution is extremely localized, the admissible adaptive scheme allows to avoid any undesired oscillation in the solution, while this feature is not maintained throughout the entire simulation when we adopt a non-admissible discretization.

%% file: sections/conclusions.tex
\section{Conclusions}\label{sec:conclusions}
We introduced a complete set of algorithms to perform admissible refinement and coarsening using THB-splines, and we successfully applied it to solve heat transfer problems with a moving heat source.
The numerical examples clearly show the advantages of the presented admissible adaptive discretization, in terms of both memory consumption and computational efficiency, with respect to uniform IGA mesh with the same level of accuracy.
Moreover, we demonstrated the importance of employing admissible meshes for these kinds of problems in order to avoid undesired oscillations which might lead to nonphysical results.
We also observed that the proposed algorithms lead to an optimal trade-off between the accuracy of the numerical results and the total number of DOFs of the system, showing better performances compared to other schemes with different refinement and coarsening algorithms, with or without taking into account grading parameters.
Finally, the error estimator, together with a robust admissible discretization, is able to capture the influence of multiple adjacent tracks on the final temperature distribution, matching one of the main requisites required by L-PBF applications.
As further outlooks we plan to extend the application of the presented discretization to 3D, non linear, and multi-physics problems in order to efficiently perform reliable simulations of AM processes.